\documentclass[12pt, twoside]{article}
\usepackage{amsmath,amsthm,amssymb}
\usepackage{times}
\usepackage{enumerate}
\usepackage{hyperref}
\usepackage{url}
\usepackage{verbatim}

\def\titlerunning#1{\gdef\titrun{#1}}
\makeatletter
\def\author#1{\gdef\autrun{\def\and{\unskip, }#1}\gdef\@author{#1}}
\def\address#1{{\def\and{\\\hspace*{18pt}}\renewcommand{\thefootnote}{}%
\footnote {#1}}%
\markboth{\autrun}{\titrun}}
\makeatother
\def\email#1{e-mail: #1}
\def\subjclass#1{{\renewcommand{\thefootnote}{}%
\footnote{\emph{Mathematics Subject Classification (2010):} #1}}}
\def\keywords#1{\par\medskip
\noindent\textbf{Keywords.} #1}

\newtheorem{thm}{Theorem}[section]

\newtheorem{conj}[thm]{Conjecture}
\newtheorem{lem}[thm]{Lemma}

\newtheorem{prop}[thm]{Proposition}

\newtheorem{mainthm}[thm]{Main Theorem}

\theoremstyle{definition}
\newtheorem{defin}[thm]{Definition}
\newtheorem{hyp}[thm]{Hypothesis}
\newtheorem{rem}[thm]{Remark}

\numberwithin{equation}{section}

\newcommand\AAA{\mathbb{A}}

\newcommand\CC{\mathbb{C}}

\newcommand\NN{\mathbb{N}}

\newcommand\PP{\mathbb{P}}
\newcommand\QQ{\mathbb{Q}}
\newcommand\RR{\mathbb{R}}

\newcommand\ZZ{\mathbb{Z}}

\newcommand\bD{{\mathbf{D}}}

\newcommand\bL{{\mathbf{L}}}

\newcommand\thalf{{\textstyle{\frac{1}{2}}}}

\newcommand\half{{{\frac{1}{2}}}}

\newcommand\fg{{\mathfrak{g}}}

\newcommand\fs{{\mathfrak{s}}}

\newcommand\ft{{\mathfrak{t}}}

\newcommand\eps{\varepsilon}

\newcommand\la{\lambda}
\newcommand\La{\Lambda}
\newcommand\ka{\kappa}

\newcommand\si{\sigma}
\newcommand\Si{\Sigma}

\newcommand\su{{\mathfrak{s}\mathfrak{u}}}

\newcommand\SO{\operatorname{SO}}

\newcommand\Ker{\operatorname{Ker}}

\newcommand\PD{\operatorname{PD}}

\newcommand\SW{SW}
\newcommand\Sym{\operatorname{Sym}}

\newcommand\even{{\mathrm{even}}}

\newcommand\id{{\mathrm{id}}}

\newcommand\spinc{\text{$\text{spin}^c$ }}
\newcommand\spinu{\text{$\text{spin}^u$ }}
\newcommand\Spinc{\text{$\text{Spin}^c$}}

\newcommand\sM{{\mathcal{M}}}

\begin{document}


\baselineskip=17pt


\titlerunning{Witten's Conjecture for many four-manifolds of simple type}

\title{Witten's Conjecture for many four-manifolds of simple type}

\author{Paul M. N. Feehan
\and 
Thomas G. Leness}

\date{November 22, 2014}

\maketitle

\address{P. M. N. Feehan: Department of Mathematics, Rutgers, The State University of New Jersey, Piscataway, NJ 08854-8019; \email{feehan@math.rutgers.edu}
\and
T. G. Leness: Department of Mathematics, Florida International University, Miami, FL 33199; \email{lenesst@fiu.edu}}



\begin{abstract}
We prove that Witten's Conjecture \cite{Witten} on the relationship between the Donaldson and Seiberg-Witten series for a four-manifold of Seiberg-Witten simple type with $b_1=0$ and odd $b_2^+\ge 3$ follows from our $\SO(3)$-monopole cobordism formula \cite{FL5} when the four-manifold has $c_1^2\ge \chi_h-3$ or is abundant.

\end{abstract}

\section{Introduction}
\label{sec:Introduction}
\subsection{Main results}
\label{subsec:MainResults}
Throughout this article, we shall assume that $X$ is a
{\em standard\/} four-manifold by which we mean that $X$ is
closed, connected, oriented, and smooth  with $b_1(X)=0$ and odd
$b^+(X)\ge 3$.  For such manifolds, we  define (by
analogy with their values when $X$ is a complex surface),
\begin{equation}
\label{eq:Definec1SquaredandHolcEulerChar}
c_1^2(X) := 2\chi+3\sigma
\quad\text{and}\quad
\chi_h(X) := \frac{1}{4}(\chi+\sigma),
\end{equation}
where $\chi$ and $\sigma$ are the Euler characteristic and signature of $X$.

For standard four-manifolds,
the Seiberg-Witten (SW) invariants \cite{MorganSWNotes},
\cite{SalamonSWBook}, \cite{Witten}  comprise a function with finite support,
$SW_X:\Spinc(X)\to\ZZ$, where $\Spinc(X)$ is the set of isomorphism classes
of \spinc structures on $X$.
The set of Seiberg-Witten (SW) basic classes, $B(X)$, is the image under a map
$c_1:\Spinc(X)\to H^2(X;\ZZ)$ of the support of $\SW_X$ \cite{Witten}.
A standard four-manifold $X$ has Seiberg-Witten simple type if
$c_1^2(\fs)=c_1^2(X)$ for all $c_1(\fs)\in B(X)$ and is
{\em abundant\/} if $B(X)^\perp\subset H^2(X;\ZZ)$ contains a hyperbolic summand,
where $B(X)^\perp$ denotes the orthogonal complement of $B(X)$ with respect to the intersection
form $Q_X$ on $H^2(X;\ZZ)$. We extend $Q_X$ from $H_2(X;\ZZ)$ to $H_2(X;\RR)$ by linearity.

We refer to \cite{KMStructure}, or \S \ref{subsec:Donaldsonseries} in this article, for the definitions of the Donaldson series, $\bD_X^w(h)$, Kronheimer-Mrowka (KM) basic classes, and four-manifolds of Kronheimer-Mrowka (KM) simple type.

\begin{conj}[Witten's Conjecture]
\label{conj:WittenSimpleType}
\cite{Witten}
Let $X$ be a standard four-manifold with Seiberg-Witten simple type. The four-manifold $X$ then has
Kronheimer-Mrowka simple type and the Kronheimer-Mrowka and Seiberg-Witten
basic classes coincide. For any $w\in H^2(X;\ZZ)$ and $h\in H_2(X;\RR)$, one has
\begin{equation}
\label{eq:WConjecture}
\bD^{w}_X(h)
=
2^{2-(\chi_h-c_1^2)}e^{Q_X(h)/2}
\sum_{\fs\in\Spinc(X)}(-1)^{\half(w^{2}+c_1(\fs)\cdot w)}
SW_X(\fs)e^{\langle c_1(\fs),h\rangle}.
\end{equation}
\end{conj}

E. Witten derived Formula \eqref{eq:WConjecture} using arguments from quantum field theory which, as far as the authors can tell, have no direct, mathematically rigorous justification. Consequently, the challenge ever since the publication of \cite{Witten} has been to provide a mathematically rigorous proof of Formula \eqref{eq:WConjecture}. 

In \cite{FL5}, we proved that a formula (restated in this article in
Theorem \ref{thm:Cobordism}) relating Donaldson and
Seiberg-Witten invariants followed from
certain properties, described in Remark \ref{rmk:GluingThmProperties}, of the gluing map
for $\SO(3)$ monopoles constructed in \cite{FL3}.
A proof of the required $\SO(3)$-monopole gluing-map properties is currently being developed by the authors.
The formula in Theorem \ref{thm:Cobordism} involves polynomials with unknown
coefficients depending on topological data and thus lacks
the elegance and simplicity of the formula in Conjecture \ref{conj:WittenSimpleType};
moreover, it appears extremely difficult, it not impossible, to compute these coefficients directly by the method of proof of Theorem \ref{thm:Cobordism}.
However, in this article, we use a
family of manifolds constructed by R. Fintushel, J. Park, and R. J. Stern
in \cite{FSParkSympOneBasic} to determine 
sufficiently many of these coefficients to
prove the

\begin{mainthm}
\label{thm:WittenSimpleType}
Let $X$ be a standard four-manifold with Seiberg-Witten simple type which is abundant or has $c_1^2(X)\ge \chi_h(X)-3$.
Then the $\SO(3)$-monopole cobordism formula (Theorem \ref{thm:Cobordism}) 
implies that Conjecture \ref{conj:WittenSimpleType} holds for $X$.
\end{mainthm}

The quantum field theory argument giving Witten's Formula
\eqref{eq:WConjecture} for standard four-manifolds has been extended by G. Moore and
E. Witten \cite{MooreWitten} to allow $b^+(X) \ge 1$, and $b_1(X)\ge 0$, and
four-manifolds $X$ of non-simple type.  The $\SO(3)$-monopole cobordism
gives a relation between the Donaldson and Seiberg-Witten invariants
for these manifolds as well and so should also
lead to a proof of  Moore and Witten's more general conjecture.
However, the methods of this article do not extend to the more general
case because of the lack of examples of four-manifolds not of simple type.

A proof of Witten's Conjecture, also assuming Theorem \ref{thm:Cobordism},
for a more restricted class of manifolds has appeared previously in \cite[Corollary 7]{KMPropertyP}.
Conjecture \ref{conj:WittenSimpleType} is known to hold, by direct calculation of both sides of Equation \eqref{eq:WConjecture}, for elliptic surfaces by work of R. Fintushel and R. J. Stern \cite{FSRationalBlowDown}.
Conjecture \ref{conj:WittenSimpleType} also holds for all simply-connected, minimal surfaces of general type.
Indeed, Theorem \ref{thm:WittenSimpleType} implies that Witten's Conjecture holds for all abundant four-manifolds
and this includes both elliptic surfaces and surfaces of general type by \cite[Corollary A.3]{FL2a}; by the discussion in \cite[\S A.2]{FL2a}, this includes all simply-connected, closed, complex surfaces with $b^+\ge 3$. In Remark \ref{rmk:LowerCoeffProb}, we explain why the arguments used in \S \ref{sec:DetermineCoeff} of our proof of Theorem \ref{thm:WittenSimpleType} do not appear, by themselves, sufficient to allow us to remove the restriction that $X$ be abundant or have $c_1^2(X)\ge \chi_h(X)-3$.

For a complex projective surface $X$, Mochizuki \cite{Mochizuki_2009} proved a formula (see Theorem 4.1 in \cite{Goettsche_Nakajima_Yoshioka_2011}) expressing the Donaldson invariants in a form similar to that given by the $\SO(3)$-monopole cobordism formula (our Theorem \ref{thm:Cobordism}), but the coefficients are given as the residues of a generating function for integrals of $\CC^*$-equivariant cohomology classes over
the product of Hilbert schemes of points on $X$. In \cite[p. 309]{Goettsche_Nakajima_Yoshioka_2011}, L. G\"ottsche, H. Nakajima, and K. Yoshioka suggest that the coefficients in Mochizuki's formula (which remain valid for a standard four-manifold) and in our $\SO(3)$-monopole cobordism formula are the same. 
They prove an explicit formula for complex projective surfaces
relating Donaldson invariants and Seiberg-Witten invariants of four-manifolds of simple
type using Nekrasov's deformed partition function for the $N=2$ SUSY gauge theory with a single fundamental matter and
from this formula deduce Witten's Conjecture.
In \cite[p. 323]{Goettsche_Nakajima_Yoshioka_2011}, they discuss the relationship between their approach, Mochizuki's formula, and our $\SO(3)$-monopole cobordism formula. See also \cite[pp. 344--347]{Goettsche_Nakajima_Yoshioka_2008} for a related discussion concerning their wall-crossing formula for the Donaldson invariants of a four-manifold with $b^+=1$.

\subsection{Outline of the article}
In \cite{FL5}, we proved that any Donaldson invariant of a
four-manifold $X$ can be
expressed as a polynomial $p_X$ in the intersection form of $X$, namely $Q_X$,
the Seiberg-Witten basic classes of $X$ and an additional cohomology class
$\La\in H^2(X;\ZZ)$ which does not appear in Equation \eqref{eq:WConjecture}.
If $X$ has SW-simple type,
then the coefficients of $p_X$
depend only on the degree of the Donaldson invariant, $\La^2$, $\chi_h(X)$,
$c_1^2(X)$, and $c_1(\fs)\cdot\La$ for an SW-basic class, $c_1(\fs)$.
We prove Theorem \ref{thm:WittenSimpleType} by using examples of manifolds known
to satisfy Conjecture \ref{conj:WittenSimpleType} to determine
sufficiently many
of these coefficients.

In \S \ref{sec:Prelim},
we review the definitions of the Donaldson series, the Seiberg-Witten
invariants, and results on the surgical operations
of blowing up 
and blowing down  which preserve
Equation \eqref{eq:WConjecture}. In \S \ref{sec:SO(3)Formula},
we summarize the background material from \cite{FL5} required to state our $\SO(3)$-monopole cobordism formula (Theorem \ref{thm:Cobordism}).
We give the proof of Theorem \ref{thm:WittenSimpleType} in
\S \ref{sec:DetermineCoeff}.

\section{Preliminaries}
\label{sec:Prelim}
We begin by reviewing the relevant
properties of the Donaldson and Seiberg-Witten invariants.

\subsection{Seiberg-Witten invariants}
\label{subsec:SWinvariants}
As stated in the introduction, the \emph{Seiberg-Witten invariants}
defined in \cite{Witten} (see also
\cite{MorganSWNotes,NicolaescuSWNotes,SalamonSWBook}),
define a map with finite support,
$$
\SW_X:\Spinc(X)\to\ZZ,
$$
where $\Spinc(X)$ denotes the set of \spinc structures on $X$.
For a \spinc structure $\fs=(W^\pm,\rho)$
where $W^\pm\to X$ are complex rank-two bundles
and $\rho$ is a Clifford multiplication map, define
$c_1:\Spinc(X)\to H^2(X;\ZZ)$ by
$c_1(\fs)=c_1(W^+)$.  For all
$\fs\in\Spinc(X)$, the cohomology class $c_1(\fs)$ is characteristic.

The invariant $\SW_X(\fs)$ is defined by the homology class
of $M_{\fs}$, the moduli space of Seiberg-Witten monopoles.
One calls $c_1(\fs)$ a {\em Seiberg-Witten (SW) basic class\/} if $\SW_X(\fs)\neq 0$.
Define
\begin{equation}
\label{eq:SWBasic}
B(X)=
\{c_1(\fs): \SW_X(\fs)\neq 0\}.
\end{equation}
If $H^2(X;\ZZ)$ has 2-torsion, then $c_1:\Spinc(X)\to H^2(X;\ZZ)$ is not injective;
moreover, the formulas in this article often involve (real) homology and cohomology, so we define
\begin{equation}
\label{eq:DefineCohomSW}
\SW_X':H^2(X;\ZZ)\to\ZZ,
\quad
K\mapsto\sum_{\fs\in c_1^{-1}(K)}\SW_X(\fs),
\end{equation}
and set $\SW_X(K)=0$ if $K$ is not characteristic.
With this definition, Witten's Formula \eqref{eq:WConjecture} is equivalent to
\begin{equation}
\label{eq:WConjCohom}
\bD^{w}_X(h)
=
2^{2-(\chi_h-c_1^2)}e^{Q_X(h)/2}
\sum_{K\in B(X)}(-1)^{\half(w^{2}+K\cdot w)}
SW'_X(K)e^{\langle K,h\rangle}.
\end{equation}
One says that a four-manifold, $X$, has {\em Seiberg-Witten (SW) simple type\/} if $\SW_X(\fs)\neq 0$ implies that $c_1^2(\fs)=c_1^2(X)$.

As discussed in \cite[\S 6.8]{MorganSWNotes},
there is an involution on $\Spinc(X)$, $\fs\mapsto\bar\fs$, with
$c_1(\bar\fs)=-c_1(\fs)$, defined essentially by taking
the complex conjugate bundles.  By \cite[Corollary 6.8.4]{MorganSWNotes}, one has
$\SW_X(\bar\fs)=(-1)^{\chi_h(X)}\SW_X(\fs)$ and so
$B(X)$ is closed under the action of $\{\pm 1\}$ on $H^2(X;\ZZ)$.

Let $\widetilde X=X\#\bar{\CC\PP}^2$ be the blow-up of $X$.
For every $n\in\ZZ$, there is a unique $\fs_n\in\Spinc(\bar{\CC\PP}^2)$
with $c_1(\fs_n)=(2n+1)e^*$, where $e^*\in H^2(\widetilde X;\ZZ)$
is the Poincar\'e dual of the exceptional curve.
By \cite[\S 4.6.2]{NicolaescuSWNotes},
there is a bijection,
$$
\Spinc(X)\times \ZZ \to \Spinc(\widetilde X),
\quad
(\fs_X,n)\mapsto \fs_X\#\fs_n,
$$
given by a connected-sum construction with
$c_1(\fs_X\#\fs_n)=c_1(\fs_X)+(2n+1)e^*$.
Versions of the following result have appeared in
\cite{FSTurkish}, \cite[Theorem 4.6.7]{NicolaescuSWNotes},
and \cite[Theorem 14.1.1]{Froyshov_2008}

\begin{thm}[Blow-up formula for Seiberg-Witten invariants]
\label{thm:FroyshovSWBlowUp}
\cite[Theorem 14.1.1]{Froyshov_2008}
Let $X$ be a standard four-manifold and let
$\widetilde X=X\#\bar{\CC\PP}^2$ be its blow-up.
Then $\widetilde X$ has SW-simple type if and only if that is true for $X$.
If $X$ has simple type, then
\begin{equation}
\label{eq:SWBasicsOfBlowUp}
B(\widetilde X)=
\{K\pm e^*: \text{$K\in B(X)$}\},
\end{equation}
and if $K\in B(X)$, then $\SW_{\widetilde X}'(K\pm e^*)=\SW_X'(K)$.
\end{thm}

\subsection{Donaldson invariants}
\label{subsec:Donaldsonseries}
\subsubsection{Definitions and the structure theorem}
We now recall the definition \cite[\S 2]{KMStructure} of the Donaldson series
for standard four-manifolds.
For any choice of $w\in H^{2}(X;\ZZ)$, the
\emph{Donaldson invariant} is a linear function,
$$
D^{w}_{X}:\AAA(X) \to \RR,
$$
where $\AAA(X)$ is the symmetric algebra,
$$
\AAA(X) = \Sym(H_{\even}(X;\RR)).
$$
For $h\in H_2(X;\RR)$ and a generator $x\in H_0(X;\ZZ)$,
we define $D_X^w(h^{\delta-2m}x^m)=0$ unless
\begin{equation}
\label{eq:DegreeParity}
    \delta\equiv -w^{2}-3\chi_h(X)\pmod{4}.
\end{equation}
If \eqref{eq:DegreeParity} holds, then $D_X^w(h^{\delta-2m}x^m)$
is defined by pairing
cohomology classes corresponding to elements of $\AAA(X)$ with the Uhlenbeck
compactification of a moduli space of anti-self-dual $\SO(3)$
connections \cite{DonPoly}, \cite{DK}, \cite{FrM}, \cite{KMStructure}.

A four-manifold has {\em Kronheimer-Mrowka (KM) simple type\/} if for all $w\in H^2(X;\ZZ)$ and
all $z\in \AAA(X)$ one has
\begin{equation}
\label{eq:KMSimpleType}
D^{w}_{X}(x^{2}z)=4D^{w}_{X}(z).
\end{equation}
The  \emph{Donaldson series} is a formal power series,
\begin{equation}
\label{eq:DefineDonaldsonSeries}
\bD^{w}_{X}(h) = D^{w}_{X}((1+\textstyle{\frac{1}{2}} x)e^{h}),
\quad h \in H_{2}(X;\RR),
\end{equation}
which determines all Donaldson invariants for  standard manifolds of
KM-simple type.
The Donaldson series of a manifold with KM-simple type
has the following description
(see also \cite[Theorems 5.9 and 5.13]{FSStructure}
for a proof by a different method):

\begin{thm}[Structure of Donaldson invariants]
\cite[Theorem 1.7 (a)]{KMStructure}
\label{thm:KMStructure}
Let $X$ be a standard four-manifold with KM-simple
type. Suppose that some Donaldson invariant of $X$ is non-zero.  Then there
is a function,
\begin{equation}
\label{eq:KMFunctionOnCohom}
\beta_X:H^2(X;\ZZ)\to \QQ,
\end{equation}
such that $\beta_X(K)\neq 0$ for at least one and at most finitely many
classes, $K$, which are integral lifts of $w_2(X)\in H^2(X;\ZZ/2\ZZ)$ (the
{\em KM-basic classes}), and for any $w\in H^2(X;\ZZ)$, one has the
following equality of analytic functions of $h\in H_2(X;\RR)$:
\begin{equation}
\label{eq:KMFormula}
\bD^{w}_X(h)
=
e^{Q_X(h)/2}
\sum_{K\in H^2(X;\ZZ)}(-1)^{(w^{2}+K\cdot w)/2}\beta_{X}(K)e^{\langle K,h\rangle}.
\end{equation}
\end{thm}

The following lemma reduces the proof of Conjecture \ref{conj:WittenSimpleType}
to proving that Equation \eqref{eq:WConjecture} holds.

\begin{lem}
Assume the hypotheses of Theorem \ref{thm:KMStructure}.
If Equation \eqref{eq:WConjecture} holds for $X$, then the KM-basic classes and SW-basic
classes coincide.
\end{lem}

\begin{proof}
The result follows by comparing Equation \eqref{eq:WConjCohom}
(which is equivalent to Equation \eqref{eq:WConjecture})
and Equation \eqref{eq:KMFormula} and by
exploiting the linear independence of the functions
$e^{r_i t}$ for different values of $r_i$.
\end{proof}

\subsubsection{Independence from $w$}
We now discuss the role of $w$.
Proofs that the condition
\eqref{eq:KMSimpleType} is
independent of $w$ appear, in varying degrees of generality, in
\cite{FroyshovFiniteType}, \cite{KMStructure}, \cite{MunozBasicNonSimple}, \cite{WWFiniteType}:

\begin{thm}
\label{thm:KMSimpleTypeIndepOfw}
\cite{KMStructure},
\cite[Theorem 2]{MunozBasicNonSimple}
Let $X$ be a standard four-manifold.  If Equation \eqref{eq:KMSimpleType} holds for one
$w\in H^2(X;\ZZ)$, then it holds for all $w$.
\end{thm}

The following proposition allows us to work with a specific $w$:

\begin{prop}
Let $X$ be a standard four-manifold of SW-simple type. If Witten's Conjecture \ref{conj:WittenSimpleType}
holds for one $w\in H^2(X;\ZZ)$, then it holds for all $w\in H^2(X;\ZZ)$.
\end{prop}

\begin{proof}
Assume that Conjecture \ref{conj:WittenSimpleType} and hence Equation \eqref{eq:WConjCohom} holds for some $w_0\in H^2(X;\ZZ)$,
\begin{equation}
\label{eq:WC1}
\begin{aligned}
{}&
e^{\half Q_X(h)}
\sum_{K\in B(X)}(-1)^{\half(w^2_0+K\cdot w_0)}\beta_{X}(K)e^{\langle K,h\rangle}
\\
{}&\quad=
2^{2-(\chi_h-c_1^2)}
e^{\half Q_X(h)}
\sum_{K\in B(X)}(-1)^{\half(w^{2}_0+K\cdot w_0)}\SW'_{X}(K)e^{\langle K,h\rangle}.
\end{aligned}
\end{equation}
We shall denote the SW-basic classes by $K_i$, for $1\leq i\leq r$, so $B(X)=\{K_1,\dots,K_s\}$.
Because $Q_X$ is indefinite, the following subset of $H_2(X;\RR)$ is non-empty:
$$
U=Q_X^{-1}(0) \setminus
\left(\bigcup_{i<j}\ (K_i-K_j)^{-1}(0)\right) \subset H_2(X;\RR).
$$
If $r_i=\langle K_i,h\rangle$ for  some fixed $h_0\in U$,
then $r_i\neq r_j$ for $i\neq j$.
Replacing $h$ by $th_0$ where $t\in\RR$ in \eqref{eq:WC1} gives
$$
\sum_{i=1}^s (-1)^{\half(w^2_0+K_i\cdot w_0)}\beta_{X}(K_i)e^{r_it}
=
2^{2-(\chi_h-c_1^2)}
\sum_{i=1}^s (-1)^{\half(w^{2}_0+K_i\cdot w_0)}\SW'_{X}(K_i)e^{r_it}.
$$
The preceding identity and linear independence of the functions
$e^{r_1 t},\dots,e^{r_s t}$ imply that
\begin{equation}
\label{eq:KMConstAreSW}
\beta_{X}(K)=2^{2-(\chi_h-c_1^2)}\SW'_X(K).
\end{equation}
Let $w$ be any other element of $H^2(X;\ZZ)$. 
Since $X$ has KM-simple type for $w_0$ (by our hypothesis that Conjecture \ref{conj:WittenSimpleType} holds for some $w_0$),
Theorem \ref{thm:KMSimpleTypeIndepOfw} implies that $X$ has KM-simple 
type for $w$.
The conclusion now follows from Equations \eqref{eq:KMFormula} and \eqref{eq:KMConstAreSW}.
\end{proof}

\subsubsection{Behavior under blow-ups}
We note that the KM-simple type condition \eqref{eq:KMSimpleType} is invariant under blow-ups.

\begin{prop}
\label{prop:KMBlowUpInvariance}
A standard four-manifold $X$ has KM-simple type if and only if
its blow-up $\widetilde X$ has KM-simple type.
\end{prop}

\begin{proof}
Assume $\widetilde X$ has KM-simple type.  The blow-up
formula $D^w_X(z)=D^w_{\tilde X}(z)$ provided by \cite[Theorem III.8.4]{FrM}
implies that,
for any $z\in\AAA(X)$,
$$
D^w_X(x^2z)=D^w_{\widetilde X}(x^2z)=4D^w_{\widetilde X}(z)=D^w_X(z),
$$
and thus $X$ has KM-simple type. The converse implication follows from \cite[Proposition 1.9]{KMStructure}.
\end{proof}

We also note the behavior of Witten's Formula \eqref{eq:WConjecture} under blow-up.

\begin{thm}
\label{thm:WCBlowDownInvariance}
\cite[Theorem 8.9]{FSRationalBlowDown}
Let $X$ be a standard four-manifold. Then Witten's Formula \eqref{eq:WConjecture} holds for $X$ if and only if it holds for the blow-up, $\widetilde X$.
\end{thm}

\subsubsection{Donaldson invariants determined by Witten's Formula}
Theorem \ref{thm:KMStructure} gives the following values for
Donaldson invariants of four-manifolds satisfying Conjecture
\ref{conj:WittenSimpleType}.
For a standard four-manifold, $X$, we define
\begin{equation}
\label{eq:Defn_c(X)}
c(X) := \chi_h(X)-c_1^2(X), 
\end{equation} 
where $\chi_h(X)$ and $c_1^2(X)$ are given in \eqref{eq:Definec1SquaredandHolcEulerChar}.
 
\begin{lem}
\label{lem:DInvarForWSTManifolds}
Let $X$ be a standard four-manifold.
Then Witten's Formula \eqref{eq:WConjecture} holds and $X$ has KM-simple type
if and only if the Donaldson invariants of $X$ satisfy
$D^w_X(h^{\delta-2m}x^m)=0$, when $\delta$ does not obey \eqref{eq:DegreeParity}, and when $\delta$ obeys \eqref{eq:DegreeParity}, then
\begin{equation}
\begin{aligned}
\label{eq:DInvarForWC}
{}&
D^w_X(h^{\delta-2m}x^m)
\\
{}&\quad
=
\sum_{\begin{subarray}{l}i+2k\\=\delta-2m\end{subarray}}
\sum_{K\in B(X)}
(-1)^{{\eps(w,K)}}
\frac{\SW'_X(K) (\delta-2m)!}{2^{k+c(X)-2-m} k!i!}
 \langle K,h\rangle^i Q_X^k(h),
\end{aligned}
\end{equation}
where $\eps(w,K):=\thalf(w^2+w\cdot K)$.
\end{lem}

\begin{proof}
Assume that Witten's Formula \eqref{eq:WConjecture}, and hence Equation \eqref{eq:WConjCohom}, holds and that $X$ has KM-simple type.
By definition, the Donaldson invariant,
$D^w_X(h^{\delta-2m}x^m)$,
will vanish unless $\delta$ obeys \eqref{eq:DegreeParity}.
Then Equation \eqref{eq:WConjCohom} holds for $X$ if and only if
\begin{align*}
{}&
2^{c(X)-2}\sum_{d=0}^\infty \frac{1}{d!}D^w_X(h^d) +  \frac{1}{d!}D^w_X(h^dx)
\\
{}&\quad=
\left( \sum_{k=0}^\infty \frac{1}{2^k k!} Q_X^k(h)\right)
\left(
\sum_{i=0}^\infty
\frac{1}{i!}
\sum_{K\in B(X)}
(-1)^{\eps(w,K)} \SW'_X(K)\langle K,h\rangle^i
\right)
\\
{}&\quad=
\sum_{d=0}^\infty
\sum_{i+2k=d}
\sum_{K\in B(X)}
(-1)^{\eps(w,K)}
\frac{\SW'_X(K)}{2^k k!i!}
 \langle K,h\rangle^i Q_X^k(h).
\end{align*}
The parity restriction \eqref{eq:DegreeParity} implies
that, for $d\not\equiv-w^2-3\chi_h\pmod 2$, one has
$$
D^w_X(h^d)+\frac{1}{2} D^w_X(h^d x)=0,
$$
while, for $d\equiv-w^2-3\chi_h\pmod 2$, 
Equation \eqref{eq:WConjCohom} holds for $X$ if and only if
\begin{align*}
{}&
2^{c(X)-2}
\left(
D^w_X(h^d)+\frac{1}{2} D^w_X(h^d x)
\right)
\\
{}&\quad=
\sum_{i+2k=d}
\sum_{K\in B(X)}
(-1)^{\eps(w,K)}
\frac{\SW'_X(K)d!}{2^k k!i!}
\langle K,h\rangle^iQ_X^k(h) .
\end{align*}
We can now read off the value of $D^w_X(h^{\delta-2m}x^m)$ from the preceding equation as follows.
If $\delta\equiv -w^2-3\chi_h\pmod 4$ and $m$ is even, then
$\delta-2m\equiv -w^2-3\chi_h\pmod 4$ so, by the KM-simple type condition \eqref{eq:KMSimpleType}
and the vanishing condition \eqref{eq:DegreeParity}
(which implies that the term $D^w_X(h^{\delta-2m}x)$ below is zero),
\begin{align*}
D^w_X(h^{\delta-2m}x^m)
{}&=2^m \left( D^w_X(h^{\delta-2m})+\frac{1}{2}  D^w_X(h^{\delta-2m}x)\right)
\\
{}&=
\sum_{\begin{subarray}{l}i+2k\\=\delta-2m\end{subarray}}
\sum_{K\in B(X)}
(-1)^{\eps(w,K)}
\frac{\SW'_X(K) (\delta-2m)!}{2^{k+c(X)-2-m} k!i!}
\langle K,h\rangle^i Q_X^k(h).
\end{align*}
Similarly, if
$\delta\equiv -w^2-3\chi_h\pmod 4$ and $m$ is odd, then
$\delta-2m+2\equiv -w^2-3\chi_h\pmod 4$ so, by the KM-simple type condition 
and the vanishing condition \eqref{eq:DegreeParity},
\begin{align*}
D^w_X(h^{\delta-2m}x^m)
{}&=2^{m-1} D^w_X(h^{\delta-2m}x)
\\
{}&=2^m \left( D^w_X(h^{\delta-2m})+\frac{1}{2}  D^w_X(h^{\delta-2m}x)\right)
\\
{}&=
\sum_{\begin{subarray}{l}i+2k\\=\delta-2m\end{subarray}}
\sum_{K\in B(X)}
(-1)^{\eps(w,K)}
\frac{\SW'_X(K) (\delta-2m)!}{2^{k+c(X)-2-m} k!i!}
 \langle K,h\rangle^i Q_X^k(h),
\end{align*}
as required.

Conversely, if the Donaldson invariants satisfy Equation \eqref{eq:DInvarForWC}
then the KM-simple type condition \eqref{eq:KMSimpleType} follows immediately.  
The fact that Witten's Formula \eqref{eq:WConjecture} holds for $X$ follows by reversing the preceding arguments.
\end{proof}

\section{The SO(3) monopole cobordism formula}
\label{sec:SO(3)Formula}
In this section, we review the $\SO(3)$-monopole cobordism formula.
More detailed expositions appear in \cite{FL5,FLGeorgia,FL2a,FL2b,FLMcMaster}.

Recall that we denote \spinc structures on $X$ by $\fs=(W^\pm,\rho)$, so
$W=W^+\oplus W^-\to X$ is a rank-four, complex Hermitian vector bundle and $\rho$
is a Clifford multiplication map.
We call $\ft=(W\otimes E,\rho\otimes\id_E)$ a \emph{\spinu structure\/} if
$(W,\rho)$ is a \spinc structure and $E\to X$ is a
rank-two complex Hermitian vector bundle.
A \spinu structure, $\ft$, defines an associated bundle, $\fg_\ft=\su(E)$,
and characteristic classes
$$
c_1(\ft)=c_1(W^+)+c_1(E)\quad\text{and}\quad p_1(\ft)=p_1(\fg_\ft).
$$
We denote
\begin{equation}
\label{eq:Defn_Lambda_kappa_w} 
\La := c_1(\ft), \quad \ka := -\frac{1}{4}\langle p_1(\ft),[X]\rangle, \quad\hbox{and}\quad w=c_1(E).
\end{equation}  
We let $\sM_{\ft}$ denote the moduli space of $\SO(3)$ monopoles for the 
\spinu structure $\ft$, as defined in \cite[Equation (2.33)]{FL2a}. 
We use the class $w$ to provide an orientation for $\sM_{\ft}$.
The moduli space $\sM_{\ft}$ admits an $S^1$ action
with fixed point subspaces
given by $M^w_\ka$, the moduli space of anti-self-dual connections
on the bundle $\fg_{\ft}$, and by  Seiberg-Witten moduli
spaces, $M_{\fs}$, where $E=L_1\oplus L_2$
and $\fs=W\otimes L_1$.  For a \spinc structure, $\fs$,
with $M_\fs\subset\sM_{\ft}$, we have $(c_1(\fs)-\La)^2=p_1(\ft)$.

The dimension of $M^w_\ka$ is given by $2\delta$, where
$$
\delta=-p_1(\ft)-3\chi_h.
$$
The dimension of $\sM_{\ft}$ is $2\delta+2n_a(\ft)$, where
$n_a(\ft)$ is the complex index of a Dirac operator defined by
$\ft$ and $n_a=( I(\La)-\delta)/4$, with
\begin{equation}
\label{eq:DefineIofLa}
I(\La)=\La^2 - \frac{1}{4}(3\chi(X)+7\si(X))
=\La^2+5\chi_h(X)-c_1^2(X).
\end{equation}
Thus,
$M^w_\ka$ has positive codimension in $\sM_{\ft}$ if and only if
$I(\La)>\delta$.  Note also that because $n_a$ is an integer,
$I(\La)\equiv \delta\pmod 4$ so, recalling that $c(X)=\chi_h(X)-c_1^2(X)$,
\begin{equation}
\label{eq:La2c(X)deltaRelation}
\La^2+c(X)\equiv \delta\pmod 4,
\end{equation}
where we used the fact that $I(\La)=\La^2+c(X)+4\chi_h(X)$ from \eqref{eq:DefineIofLa}.

The moduli space $\sM_{\ft}$ is not compact but admits a type of Uhlenbeck compactification,
$$
\bar\sM_{\ft}
\subset
\cup_{\ell=0}^N\
\sM_{\ft(\ell)}\times\Sym^\ell(X),
$$
where $\ft(\ell)$ is the \spinu structure satisfying $c_1(\ft(\ell))=c_1(\ft)$
and $p_1(\ft(\ell))=p_1(\ft)+4\ell$,
\cite[Theorem 4.20]{FL1}. The $S^1$ action extends continuously
over $\bar\sM_{\ft}$.  The closure of $M^w_\ka$ in $\bar\sM_\ft$ is the
usual Uhlenbeck compactification, $\bar M^w_\ka$, of $M^w_\ka$ \cite{DK}.
There are additional fixed points
of the $S^1$ action in $\bar\sM_{\ft}$ of the form
$M_{\fs}\times \Sym^\ell(X)$.
If $\bar\bL^w_{\ft,\ka}$ and $\bar\bL_{\ft,\fs}$ are the links of
$\bar M^w_\ka$ and $M_{\fs}\times\Sym^\ell(X)$, respectively,
in $\bar\sM_{\ft}/S^1$, then $\bar\sM_{\ft}/S^1$ defines a
compact, orientable cobordism
between  $\bar\bL^w_{\ft,\ka}$
and the union, over $\fs\in\Spinc(X)$, of the links $\bar\bL_{\ft,\fs}$ .
If $I(\La)>\delta$, then pairing certain cohomology
classes with the link $\bar\bL^w_{\ft,\ka}$ gives a multiple of the Donaldson
invariant (see \cite[Proposition 3.29]{FL2b}).
As these cohomology classes are defined on the complement of
the fixed point set in $\bar\sM_{\ft}/S^1$, the cobordism gives an
equality between this multiple of the Donaldson invariant and the
pairing of these cohomology classes with
the union, over $\fs\in\Spinc(X)$, of the links $\bar\bL_{\ft,\fs}$.
In \cite{FL5}, we computed 
an expression for this pairing, giving
a cobordism formula.

\begin{hyp}[Properties of local $\SO(3)$-monopole gluing maps]
\label{hyp:Local_gluing_map_properties}
The local gluing map, constructed in \cite{FL3}, gives a continuous parametrization of a neighborhood of $M_{\fs}\times\Si$ in $\bar\sM_{\ft}$ for each smooth stratum $\Si\subset\Sym^\ell(X)$.  
\end{hyp}

Hypothesis \ref{hyp:Local_gluing_map_properties} is recorded, in greater detail, as Conjecture 6.7.1 in \cite{FL5}. The question of how to assemble the \emph{local} gluing maps for neighborhoods of $M_{\fs}\times \Si$ in $\bar\sM_{\ft}$, as $\Si$ ranges over all smooth strata of $\Sym^\ell(X)$, into a \emph{global} gluing map for a neighborhood of $M_{\fs}\times \Sym^\ell(X)$ in $\bar\sM_{\ft}$ is itself difficult --- involving the so-called `overlap problem' described in \cite{FLMcMaster} --- but one which we do solve in \cite{FL5}. See Remark \ref{rmk:GluingThmProperties} for a further discussion of this point.

\begin{thm}[$\SO(3)$-monopole cobordism formula]
\cite{FL5}
\label{thm:Cobordism}
Let $X$ be a standard four-manifold
of Seiberg-Witten simple type. Assume that Hypothesis \ref{hyp:Local_gluing_map_properties} holds.
Assume further that $w,\La\in H^2(X;\ZZ)$ and $\delta,m\in\NN$ satisfy
\begin{enumerate}
\item
$w-\La\equiv w_2(X)\pmod 2$,
\item
$I(\La)>\delta$, where $I(\La)$
is defined in \eqref{eq:DefineIofLa},
\item
$\delta\equiv -w^2-3\chi_h\pmod 4$,
\item
$\delta-2m\ge 0$.
\end{enumerate}
Then, for any $h\in H_2(X;\RR)$ and
generator $x\in H_0(X;\ZZ)$, we have
\begin{equation}
\label{eq:MainEquation}
\begin{aligned}
{}&
D^w_X(h^{\delta-2m}x^m)
\\
{}&
\quad=
\sum_{K\in B(X)}
(-1)^{\thalf(w^2-\si)+\thalf(w^2+(w-\La)\cdot K)}SW'_X(K)
f_{\delta,m}(\chi_h,c_1^2,K,\La)(h),
\end{aligned}
\end{equation}
where
the map,
$$
f_{\delta,m}(h):\ZZ\times\ZZ\times H^2(X;\ZZ)\times H^2(X;\ZZ) \to \QQ[h],
$$
taking values in the ring of polynomials in the variable $h$ with rational coefficients, is universal (independent of $X$) and given by
\begin{equation}
\label{eq:Coefficients}
\begin{aligned}
{}&
f_{\delta,m}(\chi_h,c_1^2,K,\La)(h)
\\
{}&\quad:=
\sum_{\begin{subarray}{l}i+j+2k\\=\delta-2m\end{subarray}}
a_{i,j,k}(\chi_h,c_1^2,K\cdot\La,\La^2,m)
\langle K,h\rangle^i
\langle \La,h\rangle^j
Q_X^k(h),
\end{aligned}
\end{equation}
and,
for each triple of non-negative integers, $i, j, k \in \NN$,
the coefficients,
$$
a_{i,j,k}:\ZZ\times\ZZ\times\ZZ\times\ZZ\times\NN \to \QQ,
$$
are real analytic (independent of $X$) in the variables $\chi_h$, $c_1^2$, $c_1(\fs)\cdot\La$, $\La^2$, and $m$ with rational coefficients.
\end{thm}

\begin{rem}
\label{rmk:GluingThmProperties}
The proof of Theorem \ref{thm:Cobordism} in \cite{FL5} assumes the hypothesis
\cite[Conjecture 6.7.1]{FL5}
that the local gluing map for a neighborhood of $M_{\fs}\times \Si$ in $\bar\sM_{\ft}$ gives a continuous parametrization of a neighborhood of
$M_{\fs}\times\Si$ in $\bar\sM_{\ft}$, for each smooth stratum $\Si\subset\Sym^\ell(X)$. These local gluing maps are the analogues for
$\SO(3)$ monopoles of the local gluing maps for anti-self-dual $\SO(3)$ connections constructed by
Taubes in \cite{TauSelfDual, TauIndef, TauFrame} and Donaldson and Kronheimer in \cite[\S 7.2]{DK}; see also \cite{MorganMrowkaTube, MrowkaThesis}.
We have established the existence of local gluing maps in \cite{FL3}
and expect that a proof of the continuity for the local gluing maps with
respect to Uhlenbeck limits should be similar to our proof in \cite{FLKM1}
of this property for the local gluing maps for anti-self-dual $\SO(3)$ connections.
The remaining properties of local gluing maps assumed in \cite{FL5} are
that they are injective and also surjective in the sense that
elements of $\bar\sM_{\ft}$ sufficiently close (in the Uhlenbeck topology) to $M_{\fs}\times\Si$ are
in the image of at least one of the  local gluing maps.
In special cases, proofs of these properties for the local gluing maps for anti-self-dual $\SO(3)$ connections (namely, continuity with respect to Uhlenbeck limits, injectivity, and surjectivity) have been given in
\cite[\S 7.2.5, 7.2.6]{DK}, \cite{TauSelfDual, TauIndef, TauFrame}. The authors are currently developing a proof of the required properties for the local gluing maps for $\SO(3)$ monopoles. Our proof will also yield the analogous properties for the local gluing maps for anti-self-dual $\SO(3)$ connections.
\end{rem}

\begin{rem}
\label{rmk:PropertyP}
In \cite{KMPropertyP}, Kronheimer and Mrowka show that Theorem \ref{thm:Cobordism},
together with their work on the structure of the Donaldson invariants for manifolds of simple type \cite{KMStructure}, can be used to prove that Witten's Conjecture \ref{conj:WittenSimpleType} holds for
a suitably restricted class of standard four-manifolds \cite[Corollary 7]{KMPropertyP}
and hence prove the Property P conjecture for knots. Kronheimer and Mrowka also gave a proof of Property P which did not rely on Theorem \ref{thm:Cobordism} --- see \cite[Corollary 7.23]{KMKnotsSuturesExcisions}.
\end{rem}

\section{Determining the coefficients}
\label{sec:DetermineCoeff}
In this section, we prove that a standard four-manifold $X$
of Seiberg-Witten simple type satisfying Witten's Conjecture can determine
sufficiently many of the coefficients of the polynomial,
$$
f_{\delta,m}(\chi_h,c_1^2,c_1(\fs),\La)
$$
appearing in Equation
\eqref{eq:MainEquation} with  $\chi_h=\chi_h(X)$ and $c_1^2=c_1^2(X)$
to prove Conjecture \ref{conj:WittenSimpleType},
provided $X$ is abundant or has $c_1^2(X) \geq \chi_h(X)-3$.

\subsection{Algebraic preliminaries}
We begin with a generalization of \cite[Lemma VI.2.4]{FrM}, 
which we shall later use to determine the coefficients in Equation \eqref{eq:MainEquation}.

\begin{lem}
\label{lem:AlgCoeff}
Let $V$ be a finite-dimensional real vector space.
Let $T_1,\dots,T_n $ be linearly independent elements of the dual
space $V^*$.
Let $Q$ be a quadratic form on $V$ which is non-zero
on $\cap_{i=1}^n\Ker(T_i)$.  Then $T_1,\dots,T_n,Q$ are algebraically
independent in the sense that if
$F(z_0,\dots,z_n)\in \RR[z_0,\dots,z_n]$
and $F(Q,T_1,\dots,T_n):V\to\RR$ is the zero map, then $F(z_0,\dots,z_n)$
is the zero element of $\RR[z_0,\dots,z_n]$.
\end{lem}

\begin{proof}
We use induction on $n$.  For $n=1$, the result follows from
\cite[Lemma VI.2.4]{FrM}.

Assume that there is a polynomial $F(z_0,\dots,z_n)$
such that $F(Q,T_1,\dots,T_n):V\to\RR$
is the zero map.  Assigning $z_0$ degree two and $z_i$ degree one for $i>0$,
we can assume that $F$ is homogeneous of degree $d$.  Write
$F(z_0,\dots,z_n)=z_n^r G(z_0,\dots,z_n)$, where $z_n$ does not divide
$G(z_0,\dots,z_n)$.  Because $T_n^r G(Q,T_1,\dots,T_n)$ vanishes on $V$, 
the polynomial $G(Q,T_1,\dots,T_n)$ must vanish on the dense set $T_n^{-1}(\RR^*)$
and hence on $V$.  We now write
$G(z_0,\dots,z_n)=\sum_{i=0}^m G_i(z_0,\dots,z_{n-1})z_n^{m-i}$.
Since $z_n$ does not divide $G(z_0,\dots,z_n)$,
if $G(z_0,\dots,z_n)$ is not the zero polynomial, then
$G_m(z_0,\dots,z_{n-1})$ is not zero.  However,
as $G(Q,T_1,\dots,T_n)$ is the zero map,
the function $G_m(Q,T_1,\dots,T_{n-1})$ vanishes on $\Ker(T_n)$.
If there are 
scalars
$c_1,\dots,c_{n-1}\in\RR$ such that
the restriction of $c_1T_1+\dots c_{n-1}T_{n-1}$ to $\Ker(T_n)$
vanishes, then there is 
a scalar
$c_n\in\RR$ such that
$c_1T_1+\dots +c_{n-1}T_{n-1}=c_nT_n$.  Consequently, the linear independence
of $T_1,\dots,T_n$ implies that $c_1=\dots=c_n=0$.  Hence,
the restrictions of $T_1,\dots, T_{n-1}$
to $\Ker(T_n)$ are linearly independent.  Induction then implies
that $G_m(z_0,\dots,z_{n-1})=0$, a contradiction to
$G(z_0,\dots,z_n)$ being non-zero.  Hence, $F$ must be the zero polynomial.
\end{proof}

Being closed under the action of $\{\pm 1\}$, the set $B(X)$ is not linearly independent
over $\RR$.
Thus, in order to apply Lemma \ref{lem:AlgCoeff} to determine
the coefficients $a_{i,j,k}$ in Equation \eqref{eq:Coefficients}
from examples of manifolds satisfying Witten's Formula \eqref{eq:WConjecture}, we rewrite
the sums over $B(X)$ in Equations \eqref{eq:DInvarForWC}
and \eqref{eq:MainEquation} as sums over a smaller set of
basic classes.

Let $B'(X)$ be a fundamental domain for the action of $\{\pm 1\}$ 
on $B(X)$, so the projection map, $B'(X)\to B(X)/\{\pm 1\}$,
is a bijection. Lemma \ref{lem:DInvarForWSTManifolds} 
can then be rephrased as follows.

\begin{lem}
\label{lem:ReduceDFormToB'Sum}
Let $X$ be a standard four-manifold. 
Then Witten's Formula \eqref{eq:WConjecture} holds
and $X$ has KM-simple type if and only if the
Donaldson invariants of $X$ satisfy
$$
D^w_X(h^{\delta-2m}x^m)=0, 
$$
when $\delta\not\equiv -w^2-3\chi_h\pmod 4$,
and when $\delta \equiv -w^2-3\chi_h\pmod 4$, they satisfy
\begin{equation}
\begin{aligned}
\label{eq:DInvarForWCB'Sum}
{}&
D^w_X(h^{\delta-2m}x^m)
\\
{}&\quad
=
\sum_{\begin{subarray}{l}i+2k\\=\delta-2m\end{subarray}}
\sum_{K\in B'(X)}
(-1)^{{\eps(w,K)}}
n(K)
\frac{\SW'_X(K) (\delta-2m)!}{2^{k+c(X)-3-m} k!i!}
\langle K,h\rangle^i Q_X^k(h),
\end{aligned}
\end{equation}
where
$\eps(w,K)$ is as defined in Lemma \ref{lem:DInvarForWSTManifolds} and
\begin{equation}
\label{eq:DiracSpincFunction}
n(K)
:=
\begin{cases}
    1/2, & \text{if $K=0$,}
    \\
    1, &  \text{if $K\neq 0$.}
\end{cases}
\end{equation}
\end{lem}

\begin{proof}
We will show that Equation \eqref{eq:DInvarForWC} holds if and only
if Equation \eqref{eq:DInvarForWCB'Sum} holds and so the conclusion will follow from
Lemma \ref{lem:DInvarForWSTManifolds}.

Recall from \S \ref{subsec:SWinvariants} that $K\in B(X)$ if and only if $-K\in B(X)$.
We rewrite the sum in Equation \eqref{eq:DInvarForWC}
as a sum over $B'(X)$ by combining the $K$ and $-K$ terms
as follows.  These two terms differ only in their factors
of $(-1)^{\eps(w,K)}$, and $\SW'_X(K)$, and $\langle K,h\rangle^i$.
Because $K$ is characteristic, we see that
$$
\frac{1}{2}(w^2+w\cdot K)
-
\frac{1}{2}(w^2-w\cdot K)
\equiv
w\cdot K
\equiv
w^2\pmod 2.
$$
From \cite[Corollary 6.8.4]{MorganSWNotes}, we have
$\SW'_X(-K)=(-1)^{\chi_h}\SW'_X(K)$, so
we can combine the distinct $K$ and $-K$ terms
in Equation \eqref{eq:DInvarForWC} using the identity
\begin{equation}
\label{eq:B'SignChange}
\begin{aligned}
{}&
(-1)^{\eps(w,-K)}\SW'_X(-K)\langle -K,h\rangle^i
+
(-1)^{\eps(w,K)}\SW'_X(K)\langle K,h\rangle^i
\\
{}&\quad=
\left( (-1)^{\chi_h+w^2+i}+1 \right)
(-1)^{\eps(w,K)}\SW'_X(K)\langle K,h\rangle^i.
\end{aligned}
\end{equation}
In the sum appearing in Equation \eqref{eq:DInvarForWC},
where $i+2k=\delta-2m$, we
have $i\equiv\delta\pmod 2$.
By the parity condition \eqref{eq:DegreeParity}, we have $\delta+w^2\equiv \chi_h\pmod 4$
and so $\chi_h+w^2+i\equiv \chi_h+w^2+\delta\equiv 0\pmod 2$.
Thus, if $K\neq 0$, the $K$ and $-K$ terms will combine as in
Equation \eqref{eq:B'SignChange} to give the factor of two in Equation \eqref{eq:DInvarForWCB'Sum}.
When $K=0$, the $K$ and $-K$ terms are the same
and so we must offset this factor of two using 
the expression for $n(K)$ given in \eqref{eq:DiracSpincFunction}.
\end{proof}

We now perform a similar reduction for the sum in appearing in Equation \eqref{eq:MainEquation}.
For each triple of non-negative integers, $i,j,k \in \NN$, we define a universal polynomial map,
$$
b_{i,j,k}:\ZZ\times\ZZ\times\ZZ\times\ZZ\times\NN \to \QQ,
$$
by setting
\begin{equation}
\label{eq:DefineCombinedCoeff}
\begin{aligned}
{}&
b_{i,j,k}(\chi_h,c_1^2,K\cdot\La,\La^2,m)
\\
{}&
\quad:=
(-1)^{c(X)+i}
a_{i,j,k}(\chi_h,c_1^2,-K\cdot\La,\La^2,m)
+
a_{i,j,k}(\chi_h,c_1^2,K\cdot\La,\La^2,m),
\end{aligned}
\end{equation}
where the $a_{i,j,k}$ are the
universal, rational
coefficients appearing
in the expression \eqref{eq:Coefficients}.
Definition \eqref{eq:DefineCombinedCoeff} implies that
\begin{equation}
\label{eq:CombinedCoeffReverseLaKEquality}
b_{i,j,k}(\chi_h,c_1^2,-K\cdot\La,\La^2,m)
=
(-1)^{c(X)+i}
b_{i,j,k}(\chi_h,c_1^2,K\cdot\La,\La^2,m).
\end{equation}
We also define,
\begin{equation}
\label{eq:DefineTildeEps}
\tilde\eps(w,\La,K):=\frac{1}{2}(w^2-\si)+\frac{1}{2}(w^2+(w-\La)\cdot K).
\end{equation}
We can now state the desired reduction.

\begin{lem}
\label{lem:ReduceCobordismFormToB'Sum}
Assume the hypotheses of Theorem \ref{thm:Cobordism}.
Denote the coefficients in \eqref{eq:DefineCombinedCoeff} more concisely by
$$
b_{i,j,k}(K\cdot\La) := b_{i,j,k}(\chi_h,c_1^2,K\cdot\La,\La^2,m).
$$
Then,
\begin{equation}
\begin{aligned}
\label{eq:CompareCoeff2}
D^w_X(h^{\delta-2m}x^m)
&=
\sum_{\begin{subarray}{l}i+j+2k\\=\delta-2m\end{subarray}}
\sum_{K\in B'(X)}
n(K)(-1)^{{\tilde\eps(w,\La,K)}}SW'_X(K)
\\
&\qquad\times
b_{i,j,k}(K\cdot\La)
\langle K,h\rangle^i
\langle \La,h\rangle^j
Q_X^k(h),
\end{aligned}
\end{equation}
where $n(K)$ is defined by \eqref{eq:DiracSpincFunction}.
\end{lem}

\begin{proof}
Because the class $w-\La$ is characteristic and as $K^2=c_1^2(X)$, we have
$$
\tilde\eps(w,\La,-K)
=
\tilde\eps(w,\La,K)-(w-\La)\cdot K)
\equiv
\tilde\eps(w,\La,K)+c_1^2(X)\pmod 2.
$$
For $K\neq 0$, we can combine the
distinct
$K$ and $-K$ terms
in the sum appearing in Equation \eqref{eq:MainEquation}
as in the identity \eqref{eq:B'SignChange} to
obtain the expression \eqref{eq:DefineCombinedCoeff} for the coefficients, $b_{i,j,k}$.
For $K=0$, the factor of $n(K)=1/2$ is necessary because
the addition of the two identical terms in \eqref{eq:DefineCombinedCoeff}
would correspond to counting the term for $K=-K=0$ in Equation \eqref{eq:MainEquation}
twice.
\end{proof}

\subsection{The example manifolds}
A four-manifold with the properties described in Definition \ref{defn:Useful} can be
used with Lemmas \ref{lem:AlgCoeff}, \ref{lem:ReduceDFormToB'Sum},
and \ref{lem:ReduceCobordismFormToB'Sum} to determine many of the coefficients $b_{i,j,k}$ in Equation \eqref{eq:CompareCoeff2}.

\begin{defin}[Useful four-manifolds]
\label{defn:Useful}
We call a standard four-manifold, $X$, {\em useful\/} if
\begin{enumerate}
\item
$X$ has SW-simple type, and $|B'(X)|=1$,
\item
$X$ satisfies Witten's Equation \eqref{eq:DInvarForWCB'Sum},
\item
There are cohomology classes, $f_1, f_2\in B(X)^\perp$, with
$f_i^2=0$ and $f_1\cdot f_2=1$ such that
$\{f_1,f_2\}\cup B'(X)$ is linearly independent
over $\RR$, and
\item
If $f_1,f_2$ are the cohomology classes in the previous condition,
then the restriction of $Q_X$ to
$\left(\cap_{i=1}^2\Ker(f_i)\right)\cap \left(\cap_{K\in B'(X)}\Ker(K)\right)$
is non-zero.
\end{enumerate}
\end{defin}

We prove the existence of a family of useful four-manifolds in the following lemma.

\begin{lem}[Existence of useful four-manifolds]
\label{lem:UsefulWithc=3}
For every integer $h=2,3,4,\dots,$ there is a useful four-manifold $Y_h$
with $\chi_h(Y_h)=h$, $c_1^2(Y_h)=h-3$, and $c(Y_h)=3$.
\end{lem}

\begin{proof}
In \cite[Proposition 3.5]{FSParkSympOneBasic}, R. Fintushel, J. Park, and R. Stern
construct examples of standard four-manifolds $X_p$ and $X_p'$
for integer $p\ge 4$ with $c_1^2(X_p)=2p-7$ and
$c_1^2(X_p')=2p-8$ and both satisfying $c_1^2=\chi_h-3$.
In addition, $|B(X_p)/\{\pm 1\}|=|B(X_p')/\{\pm 1\}|=1$.
The four-manifolds constructed in \cite{FSParkSympOneBasic} define a ray in the
$(\chi_h,c_1^2)$ plane but the restrictions 
on $p$ mean that
they do not include the point $\chi_h=2$ and $c_1^2=-1$.
We will write $Y_h$ for the member of this family of manifolds with $\chi_h(Y_h)=h$
and set $Y_2:=K3\#\bar{\CC\PP}^2$,  
where `$K3$' denotes the K3 surface.
We further note that $Y_3=E(3)$
by the construction in \cite[\S 3]{FSRationalBlowDown} where one notes that
the operation of rationally blowing down the empty configuration $C_1$ is trivial,
\cite{FintushelPersonalCommuncation}.
Because $B(K3)=\{0\}$ by \cite{FSRationalBlowDown},
the blow-up formula in Theorem \ref{thm:FroyshovSWBlowUp} implies that $|B'(Y_2)|=1$.

As shown 
in the discussion following Lemma 3.4 in \cite{FSParkSympOneBasic},
for $p\ge 4$, the four-manifolds $X_p$ and $X_p'$ are rational blow-downs of the elliptic surfaces
$E(2p-4)$ and $E(2p-5)$, respectively, along taut configurations
(in the sense of \cite[\S 7]{FSRationalBlowDown}) of embedded spheres.
These elliptic surfaces have SW-simple type and satisfy Conjecture
\ref{conj:WittenSimpleType} (see, for example, \cite[Theorem 8.7]{FSRationalBlowDown}).
By \cite[Theorem 8.9]{FSRationalBlowDown},  
these properties (having SW-simple type and satisfying Conjecture \ref{conj:WittenSimpleType})
are preserved under rational blowdown and hence also hold for $Y_h$ for $h>2$.
$p\ge 4$.  For $Y_2=K3\#\bar{\CC\PP}^2$, these two properties hold
because they hold for $K3=E(2)$, by \cite{KMStructure}
and \cite{FSRationalBlowDown}, and
because these properties are
preserved under blow-ups by Theorem \ref{thm:WCBlowDownInvariance}.

Recall that a four-manifold $X$ is abundant if there are cohomology classes 
$f_1,f_2\in B(X)^\perp\subset H^2(X;\ZZ)$
with $f_i^2=0$ and $f_1\cdot f_2=1$.
By \cite[Corollary A.3]{FL2a}, if $X$ is simply connected and the SW-basic classes are
all multiples of a single cohomology class, then $X$ is abundant.
This result, together with
the fact that $|B'(Y_h)|=1$  for all $h\ge 2$ implies that our four-manifolds,
$Y_h$, are abundant.

We now show that the cohomology-class linear independence property holds for the four-manifolds $Y_h$.
If the cohomology classes $f_1,f_2\in B(Y_h)^\perp$ are as described in the Definition \ref{defn:Useful} of a
useful four-manifold and $K\in B(Y_h)$ and $af_1+bf_2+cK=0$ for
some $a,b,c\in\RR$, then
$$
a=f_2\cdot(af_1+bf_2+cK)=0 \quad\hbox{and}\quad b=f_1\cdot(af_1+bf_2+cK)=0,
$$
and thus $cK=0$.  If $K\neq 0$, then $c=0$ and the set $\{K,f_1,f_2\}$ is linearly
independent.  If $K=0$, then because the four-manifolds $Y_h$ have SW-simple type,
we would have
$0=K^2=c_1^2(Y_h)$ which is only true if $h=3$ and $Y_3=E(3)$.
For $h=3$, we have $B'(Y_3)=\{F\}$, where $F$ is the Poincar\'e dual
of a generic fiber of the elliptic
fibration on $Y_3$ by \cite{FSRationalBlowDown} and $F\neq 0$.
Hence, $K\neq 0$ for all our manifolds $Y_h$, so the set
$\{K,f_1,f_2\}$ is linearly independent over $\RR$.

To prove that our manifolds $Y_h$ satisfy
the fourth condition in the Definition \ref{defn:Useful} of a useful four-manifold,
we identify the kernels of the cohomology classes $K$, $f_1$, and $f_2$
with their orthogonal complements in $H^2(Y_h;\ZZ)$ by Poincar\'e duality,
and show that the restriction of $Q_{Y_h}$ to this orthogonal complement is non-zero.
If $K^2\neq 0$, then the determinant of the restriction of $Q_{Y_h}$ to
the span of  $\{K,f_1,f_2\}$ is non-zero.  Hence, the determinant
of the restriction of $Q_{Y_h}$ (and thus the restriction of $Q_{Y_h}$)
to the orthogonal complement of this span is also non-zero.
As in the preceding paragraph, if $K^2=0$, then $h=3$ and
$Y_3=E(3)$.  If $F\in H^2(E(3);\ZZ)$ is the Poincar\'e dual of a generic
fiber of the elliptic fibration and $\si\in H^2(E(3);\ZZ)$ is
the Poincar\'e dual of a section, then
$1=F\cdot \si\equiv \si^2\pmod 2$ so
$Q_{E(3)}$ is odd and there is an isomorphism of quadratic forms
$$
\left( H^2(E(3);\ZZ), Q_{E(3)}\right)
\simeq
\left( \mathop{\oplus}\displaylimits_{i=1}^5 \ZZ e_i\right)
\oplus
\left( \mathop{\oplus}\displaylimits_{j=1}^{29} \ZZ g_j\right),
$$
where $e_i^2=1$ and $g_j^2=-1$. 
Following the argument of \cite[Lemma A.4]{FL2a},
we define
\begin{gather*}
L:=3e_1+3e_2+3e_3+e_4+e_5+\sum_{j=1}^{29} g_j,
\\
\tilde f_1:=e_5+g_2,
\quad
\tilde f_2:=e_5+g_3,
\quad
P:= e_1-e_2.
\end{gather*}
Then, $L$ is primitive and characteristic with $L^2=0$,
while $\{\tilde f_1,\tilde f_2\}$ span a hyperbolic summand orthogonal
to $L$.  The class $P$ is orthogonal to the
the span of $\{L,\tilde f_1,\tilde f_2\}$ and $P^2\neq 0$.
Thus, $Q_{E(3)}$ is non-zero on the orthogonal complement
of the span $\{L,\tilde f_1,\tilde f_2\}$.
Because $\si\cdot F=1$, then $F$ is primitive as
well as characteristic with $F^2=0$.
As observed in \cite[Lemma A.4]{FL2a},
a result of Wall (see \cite[Proposition 1.2.28]{WallUnimodQuadForms})
implies that the orthogonal group of 
$(H^2(E(3);\ZZ),Q_{E(3)})$ acts transitively on the primitive
characteristic elements with a given square.
Hence, there is an isometry of $(H^2(E(3);\ZZ),Q_{E(3)})$
mapping $L$ to $F$.  If we take $f_i$ to be the image
of $\tilde f_i$ under this isometry, then we see
that $Q_{E(3)}$ is non-zero on the orthogonal complement
of the span of $\{F, f_1, f_2\}$, as desired.
\end{proof}

\subsection{The blow-up formulas}
To determine the coefficients $b_{i,j,k}$ for a sufficiently wide
range of values of $\chi_h$, $c_1^2$, $\La^2$, and
$K\cdot\La$, we will need to
work with the blow-ups of the useful four-manifolds described in Lemma \ref{lem:UsefulWithc=3}.  Thus,
let $\widetilde X(n)$ be the blow-up of $X$ at $n$ points, where $X$ is
one of the useful four-manifolds described in Lemma \ref{lem:UsefulWithc=3}.
For non-negative integers $m\le n$, we will consider
$H^2(\widetilde X(m);\ZZ)$ as a subspace of
$H^2(\widetilde X(n);\ZZ)$ using the inclusion defined by
the pullback of the blowdown map.
Let $\{e_1,\ldots,e_n\} \subset H_2(\widetilde X(n);\ZZ)$ be the homology classes of the
exceptional curves and let $e_u^*:=\PD[e_u]$,
for $u=1,\ldots,n$.

We now describe $B(\widetilde X(n))$ in more detail.
Let $\pi_u:(\ZZ/2\ZZ)^n\to \ZZ/2\ZZ$ be projection onto the $u$-th factor.
For $K\in B(X)$ and $\varphi\in (\ZZ/2\ZZ)^n$, define
\begin{equation}
\label{eq:Kvarphi_Kzero} 
K_\varphi:=K+\sum_{u=1}^n (-1)^{\pi_u(\varphi)}e_u^* \quad\hbox{and}\quad K_0:=K+\sum_{u=1}^n e_u^*.
\end{equation} 
If $0\notin B(X)$, then the Seiberg-Witten blow-up formula \eqref{eq:SWBasicsOfBlowUp} implies that
$$
B'(\widetilde X(n))=
\{K_\varphi: K\in B'(X) \text{ and } \varphi\in (\ZZ/2\ZZ)^n\}.
$$
Even if the set $B'(X)$ of SW-basic classes is linearly independent, the set $B'(\widetilde X(n))$ will not
be linearly independent for $n\ge 2$.  

To rewrite Lemma \ref{lem:ReduceCobordismFormToB'Sum} in terms of linearly independent
SW-basic
classes, we will require a result from combinatorics.
For a function $f:\ZZ\to\RR$ and $p,q\in\ZZ$, define
\begin{equation}
\label{eq:Nabla_f_on_integers} 
(\nabla^q_pf)(x):=f(x)+(-1)^q f(x+p), \quad\forall\, x \in \ZZ,
\end{equation} 
and for $a\in\ZZ/2\ZZ$ and $p\in \ZZ$, define 
\begin{equation}
\label{eq:pa}
pa:=-\frac{1}{2}(-1 +(-1)^a)p.
\end{equation} 
We then have

\begin{lem}
\label{lem:PermutationSumAsDifferenceOperator}
Let $f:\ZZ\to\RR$ be a function and $n\geq 1$ an integer. Then, for all $(p_1,\dots,p_n)$ and $(q_1,\dots,q_n)$ in $\ZZ^n$, one has
$$
\sum_{\varphi\in(\ZZ/2\ZZ)^n}
(-1)^{\sum_{u=1}^n q_u\pi_u(\varphi)} f\left(x+\sum_{u=1}^n p_u \pi_u(\varphi)\right)
=
(\nabla^{q_1}_{p_1}\nabla^{q_2}_{p_2}\dots\nabla^{q_n}_{p_n}f)(x),
$$
and if $C$ is the constant function, then
\begin{equation}
\label{eq:ConstantDifference}
(\nabla^{q_n}_{p_n}\nabla^{q_{n-1}}_{p_{n-1}}\dots\nabla^{q_1}_{p_1}C)
=
\begin{cases}
0, & \text{if $q_u\equiv 1\pmod 2$ for some $u \in \{1,\ldots,n\}$,}
\\
2^{n}C, & \text{if $q_u\equiv 0\pmod 2$ for all $u \in \{1,\ldots,n\}$.}
\end{cases}
\end{equation}
\end{lem}

\begin{proof}
The proof uses  induction on $n$.
For $n=1$, the statement is trivial.
Define,
$$
(L^{q_1,\dots,q_n}_{p_1,\dots,p_n}f)(x)
:=
\sum_{\varphi\in(\ZZ/2\ZZ)^n}
(-1)^{\sum_{u=1}^n q_u\pi_u(\varphi)} f(x+\sum_{u=1}^n p_u \pi_u(\varphi)).
$$
For $n\ge 2$, the preceding expression can be expanded as
\begin{align*}
{}&
\sum_{\varphi\in\pi_n^{-1}(0)}
(-1)^{\sum_{u=1}^{n-1} q_u\pi_u(\varphi)} f(x+\sum_{u=1}^{n-1} p_u \pi_u(\varphi))
\\
{}&\quad\quad
+(-1)^{q_n}
\sum_{\varphi\in\pi_n^{-1}(1)}
(-1)^{\sum_{u=1}^{n-1} q_u\pi_u(\varphi)} f(x+p_n+\sum_{u=1}^{n-1} p_u \pi_u(\varphi))
\\
{}&\quad=
(L^{q_1,\dots,q_{n-1}}_{p_1,\dots,p_{n-1}}f)(x)
+(-1)^{q_n}
(L^{q_1,\dots,q_{n-1}}_{p_1,\dots,p_{n-1}}f)(x+p_n)
\\
{}&\quad=
(\nabla^{q_n}_{p_n}(L^{q_1,\dots,q_{n-1}}_{p_1,\dots,p_{n-1}}f))(x),
\end{align*}
where in the penultimate step we have identified $(\ZZ/2\ZZ)^{n-1}$
with $\pi_n^{-1}(0)$ and $\pi_n^{-1}(1)$ as sets.  
The first assertion in the lemma now follows by induction.

The identity \eqref{eq:ConstantDifference} follows from the fact that
$$
\nabla^q_p C =
C+(-1)^qC=
\begin{cases}
0, & \text{if $q\equiv 1\pmod 2$,} \\
2C, & \text{if $q\equiv 0\pmod 2$,}
\end{cases}
$$
and induction on $n$.
\end{proof}

If $X$ is a four-manifold with
blow-up $\widetilde X(n)$ for some integer $n\geq 1$ and $w\in H^2(X;\ZZ)\subset H^2(\widetilde X(n);\ZZ)$, we denote
\begin{equation}
\label{eq:Defn_tilde_w} 
\tilde w:=w+\sum_{u=1}^n w_ue_u^*.
\end{equation} 
We can now rewrite Lemmas \ref{lem:ReduceDFormToB'Sum} and
\ref{lem:ReduceCobordismFormToB'Sum} in terms of linearly independent
SW-basic
classes.

\begin{lem}
\label{lem:BlowUpCobordism}
Continue the notation of the preceding paragraphs
and Definition \ref{defn:Useful}.
Let $X$ be a useful four-manifold and $n\geq 1$ an integer.
For $w\in H^2(X;\ZZ)\subset H^2(\widetilde X(n);\ZZ)$,
let $\tilde w$ be as in \eqref{eq:Defn_tilde_w}.
Let $\La\in H^2(\widetilde X(n);\ZZ)$ satisfy
$I(\La)>\delta$
and $\La-\tilde w\equiv w_2(\widetilde X(n))\pmod 2$.
Define
$$
b_{i,j,k}( K_\varphi\cdot \La)
:=
b_{i,j,k}(\chi_h(\widetilde X(n)),c_1^2(\widetilde X(n)),K_\varphi\cdot\La,\La^2,m).
$$
Then, for $\delta-2m\ge 0$,
\begin{equation}
\label{eq:BlownUpUsefulCobordismFormula}
\begin{aligned}
(-1)^{\eps(\tilde w,K_0)}
&\sum_{\begin{subarray}{c}i_0+\dots+i_n+2k\\=\delta-2m\end{subarray}}
\binom{i_0+\dots + i_n}{i_0,\, i_1,\, \ldots, i_n}
\frac{\SW'_X(K) (\delta-2m)!}{2^{k+c+n-3-m} k!i!}
\\
&\qquad\times
p^{\tilde w}(i_1,\,i_2,\,\ldots,i_n)
\langle K,h\rangle^{i_0}
\prod_{u=1}^n\langle e_u^*,h\rangle^{i_u}
Q_X^k(h)
\\
=
(-1)^{\tilde\eps(\tilde w,\La,K_0)}
&\sum_{\begin{subarray}{c}i_0+\dots+i_n+j+2k\\=\delta-2m\end{subarray}}
\binom{i_0+\dots + i_n}{i_0,\,i_1,\, \ldots, i_n}\SW'_X(K)
\\
&\qquad\times
\sum_{\varphi\in (\ZZ/2\ZZ)^n}
(-1)^{\sum_{u=1}^n (1+i_u)\pi_u(\varphi)}
b_{i,j,k}(K_\varphi\cdot\La )
\\
{}&\qquad\times
\langle K,h\rangle^{i_0}\prod_{u=1}^n\langle e_u^*,h\rangle^{i_u}
\langle \La,h\rangle^j Q_X^k(h),
\end{aligned}
\end{equation}
where $c=c(X)=\chi_h(X)-c_1^2(X)$, 
as in \eqref{eq:Defn_c(X)}, and
\begin{equation}
\label{eq:DefineSumFactor1}
p^{\tilde w}(i_1,i_2,\dots,i_n)
:=
\begin{cases}
0, & \text{if $w_q+i_q\equiv 1\pmod 2$ for some $q \in \{1,\ldots n\}$,}
\\
2^{n}, & \text{if $w_q+i_q\equiv 0\pmod 2$ for all $q \in \{1,\ldots, q\}$.}
\end{cases}
\end{equation}
\end{lem}

\begin{proof}
Comparing Equations \eqref{eq:DInvarForWCB'Sum} and \eqref{eq:CompareCoeff2}
yields, for $\eps(\tilde w,\varphi)=\thalf(\tilde w^2+\tilde w\cdot K_\varphi)$,
\begin{equation}
\label{eq:CoeffDet1}
\begin{aligned}
{}&
\sum_{\begin{subarray}{l}i+2k\\=\delta-2m\end{subarray}}
\sum_{\varphi\in (\ZZ/2\ZZ)^n}
(-1)^{{\eps(\tilde w,\varphi)}}
\frac{\SW'_X(K) (\delta-2m)!}{2^{k+c+n-3-m} k!i!}
\langle K_\varphi,h\rangle^iQ_X^k(h)
\\
{}&\quad
=
\sum_{\begin{subarray}{l}i+j+2k\\=\delta-2m\end{subarray}}
\sum_{\varphi\in (\ZZ/2\ZZ)^n}
(-1)^{{\tilde\eps(w,\La,K_\varphi)}}SW'_X(K)
b_{i,j,k}( K_\varphi\cdot\La)
\langle K_\varphi,h\rangle^i
\langle \La,h\rangle^j
Q_X^k(h).
\end{aligned}
\end{equation}
For $\varphi\in(\ZZ/2\ZZ)^n$, we have
\begin{equation}
\label{eq:OrientationVariation}
\eps(\tilde w,\varphi)
\equiv
\frac{1}{2}
(\tilde w^2 + \tilde w\cdot K_\varphi)
\equiv
\frac{1}{2}
(\tilde w^2 + \tilde w\cdot K_0)
+
\sum_{u=1}^n w_u\pi_u(\varphi)
\pmod 2.
\end{equation}
By the multinomial theorem, for $\varphi\in(\ZZ/2\ZZ)^n$
we can expand the factor $\langle K_{\varphi},h\rangle^i$ as
\begin{equation}
\label{eq:ExpandPowerOfKphi}
\begin{aligned}
\langle K_\varphi,h\rangle^i
=
\sum_{i_0+\dots+i_n=i}
\binom{i}{i_0,\ i_1, \ \dots, i_n}
(-1)^{\sum_{u=1}^n \pi_u(\varphi)i_u}
\langle K,h\rangle^{i_0}
\prod_{u=1}^n\langle e_u^*,h\rangle^{i_u},
\end{aligned}
\end{equation}
where, for $i=i_0+\dots +i_n$,
$$
\binom{i}{i_0,\, i_1, \, \ldots, i_n}=\frac{i!}{i_0!i_1!\dots i_n!}.
$$
The identities \eqref{eq:OrientationVariation}
and \eqref{eq:ExpandPowerOfKphi} imply that we can rewrite
the left-hand side of Equation \eqref{eq:CoeffDet1} as
\begin{equation}
\label{eq:CoeffDet2a}
\begin{aligned}
{}&
\sum_{\begin{subarray}{l}i+2k\\=\delta-2m\end{subarray}}
\sum_{\varphi\in (\ZZ/2\ZZ)^n}
(-1)^{{\eps(\tilde w,\varphi)}}
\frac{\SW'_X(K)(\delta-2m)!}{2^{k+c+n-3-m} k!i!}
\langle K_\varphi,h\rangle^iQ_X^k(h)
\\
{}&\qquad=
(-1)^{\eps(\tilde w,K_0)}
\sum_{\begin{subarray}{c}i_0+\dots+i_n+2k\\=\delta-2m\end{subarray}}
\binom{i_0+\dots+i_n}{i_0, \, \ldots, i_n}
\frac{\SW'_X(K) (\delta-2m)!}{2^{k+c+n-3-m} k!i!}
\\
{}&\qquad\qquad
\times
\sum_{\varphi\in (\ZZ/2\ZZ)^n}
(-1)^{\sum_{u=1}^n \pi_u(\varphi)(w_u+i_u)}
\langle K,h\rangle^{i_0}
\prod_{u=1}^n\langle e_u^*,h\rangle^{i_u}
Q_X^k(h).
\end{aligned}
\end{equation}
By applying Lemma \ref{lem:PermutationSumAsDifferenceOperator},
we write the sum over $\varphi\in(\ZZ/2\ZZ)^n$ in Equation \eqref{eq:CoeffDet2a}
as
$$
\sum_{\varphi\in (\ZZ/2\ZZ)^n}(-1)^{\sum_{u=1}^n \pi_u(\varphi)(w_u+i_u)}
=
\nabla^{w_1+i_1}_0\cdots \nabla^{w_n+i_n}_0
1.
$$
Equation \eqref{eq:ConstantDifference} shows that the preceding expression is
equal to $p^{\tilde w}(i_1,\dots,i_n)$, as defined in \eqref{eq:DefineSumFactor1}.
Therefore, Equation \eqref{eq:CoeffDet2a}
implies that the left-hand side of Equation \eqref{eq:CoeffDet1}
equals the left-hand side of Equation \eqref{eq:BlownUpUsefulCobordismFormula}.

We now rewrite the right-hand side of Equation \eqref{eq:CoeffDet1}.
The discussion is essentially the same as that for the left-hand side.
However, note that
$$
\tilde \eps(\tilde w,\La,K_\varphi)
-
\tilde \eps(\tilde w,\La,K_0)
=
\frac{1}{2} (\La-\tilde w)\cdot(K_\varphi-K_0).
$$
Because
$$
K_\varphi-K_0
=
\sum_{u=1}^n
((-1)^{\pi_u(\varphi)}-1)e_u^*
=
-2\sum_{u=1}^n \pi_u(\varphi)e_u^*,
$$
and since $\La-\tilde w$ is characteristic, we have
$$
\frac{1}{2} (\La-\tilde w)\cdot(K_\varphi-K_0)
\equiv
\sum_{u=1}^n \pi_u(\varphi)\pmod 2.
$$
The preceding identity replaces the orientation sign-change factor
computed in \eqref{eq:OrientationVariation},
and we can conclude that the right-hand side of Equation \eqref{eq:CoeffDet1}
is equal to the right-hand side of Equation \eqref{eq:BlownUpUsefulCobordismFormula}.
\end{proof}

\subsection{Determining the coefficients $b_{i,j,k}$}
We now apply Lemmas \ref{lem:AlgCoeff} and \ref{lem:BlowUpCobordism}
to the manifolds discussed in Lemma \ref{lem:UsefulWithc=3} to
determine the coefficients $b_{i,j,k}$ with $i\ge c(X)-3>0$.

\begin{prop}
\label{prop:HighDegreeCoefficients}
For any integers $x,y$ and for any integers $m\ge 0$, $n>0$, and $\chi_h\ge 2$
and for any non-negative integers $i$, $j$, $k$ satisfying
$i+j+2k=\delta-2m$, $i\ge n$, and $2y>\delta-4\chi_h-3-n$,
the coefficients
$b_{i,j,k}(\chi_h,c_1^2,K\cdot\La,\La^2,m)$
defined in Equation \eqref{eq:DefineCombinedCoeff} satisfy
$$
b_{i,j,k}(\chi_h,\chi_h -3-n,2x,2y,m)
=
\begin{cases}
(-1)^{x+y}
\displaystyle\frac{(\delta-2m)!}{k!i!}2^{m-k-n},  & \text{if $j=0$,}
\\
0, & \text{if $j>0$.}
\end{cases}
$$
\end{prop}

\begin{proof}
For one of the useful four-manifolds,  $X$, described in
Lemma \ref{lem:UsefulWithc=3}, let $\widetilde X(n)$ be the blow-up of $X$ at $n$ points.
We apply Lemma \ref{lem:BlowUpCobordism} with
$$
\La=(y+2x^2)f_1+f_2+2x e_1^*,
$$
where $f_1,f_2\in B(X)^\perp$
are the cohomology classes in Definition \ref{defn:Useful}
satisfying $f_i^2=0$ and $f_1\cdot f_2=1$.
Thus,
$$
\La^2=2y
\quad\hbox{and}\quad
K_0\cdot\La=-2x.
$$
The condition $2y>\delta-4\chi_h-3-n$ implies that $I(\La)>\delta$.
Observe that
$$
(K_\varphi-K_0)\cdot \La
=
\begin{cases}
0, &\text{if $\pi_1(\varphi)=0$,} \\
4x, &\text{if $\pi_1(\varphi)=1$.}
\end{cases}
$$
If we write $\tilde w=w+\sum_{u=1}^n w_ue_u^*$, as in \eqref{eq:Defn_tilde_w}, then the requirement
that $\La-\tilde w$ is characteristic implies that $w_u\equiv 1\pmod 2$
for all $u$.  Hence, the coefficient of the term,
\begin{equation}
\label{eq:GeneralTerm}
K^{i_0}(e_1^*)^{i_1}\cdots (e_n^*)^{i_n}\La^j Q_X^k,
\end{equation}
on the left-hand side of Equation \eqref{eq:BlownUpUsefulCobordismFormula}
will vanish if $j>0$ 
while, if $j=0$, the coefficient is equal to
\begin{equation}
\label{eq:LHSCoefficientForBigi}
(-1)^{\eps(\tilde w,K_0)}
\binom{i}{i_0,i_1,\dots,i_n}
\frac{\SW'_X(K) (\delta-2m)!}{2^{k+n-m} k!i!}
p^{\tilde w}(i_i,\dots,i_n),
\end{equation}
where $i=i_0+\dots + i_n$.

The coefficient of the term \eqref{eq:GeneralTerm}
on the right-hand side of Equation \eqref{eq:BlownUpUsefulCobordismFormula} is
\begin{equation}
\label{eq:RHSCoeff2a}
\begin{aligned}
{}&
(-1)^{\tilde \eps(\tilde w,\La,K_0)}
\binom{i}{i_0,\,\ldots,i_n}\SW'_X(K)
\\
&\quad\ \
\times\left(
b_{i,j,k}(-2x)
\left(
\sum_{\varphi\in\pi_1^{-1}(0)} (-1)^{\sum_{u=1}^n(1+i_u)\pi_u(\varphi)}
\right)\right.
\\
&\quad\qquad
+
b_{i,j,k}(2x)
\left.\left(
\sum_{\varphi\in\pi_1^{-1}(1)} (-1)^{\sum_{u=1}^n(1+i_u)\pi_u(\varphi)}
\right)
\right).
\end{aligned}
\end{equation}
Equation \eqref{eq:ConstantDifference} implies that, for $a=0,1$,
\begin{align*}
\sum_{\varphi\in \pi^{-1}(a)}(-1)^{\sum_{u=2}^n(1+i_u)\pi_u(\varphi)}
{}&=
\begin{cases}
0, & \text{if $i_q\equiv 0\pmod 2$ for some $q\in\{2,\ldots,n\}$,}
\\
2^{n-1}, &\text{if $i_q\equiv 1\pmod 2$ for all $q\in\{2,\ldots,n\}$.}
\end{cases}
\end{align*}
We define a map $p^1:\ZZ^{n-1}\to\ZZ$ by setting $p^1(i_2,\dots,i_n)$ equal to the right-hand side of the preceding expression.
Hence,
\begin{align*}
\sum_{\varphi\in\pi_1^{-1}(0)} (-1)^{\sum_{u=1}^n(1+i_u)\pi_u(\varphi)}
{}&=
p^1(i_2,\dots,i_n),
\\
\sum_{\varphi\in\pi_1^{-1}(1)} (-1)^{\sum_{u=1}^n(1+i_u)\pi_u(\varphi)}
{}&=
(-1)^{1+i_1}p^1(i_2,\dots,i_n).
\end{align*}
The identity \eqref{eq:CombinedCoeffReverseLaKEquality} and
the identity $\La^2-\delta\equiv c(\widetilde X(n))\pmod 4$
implied by \eqref{eq:La2c(X)deltaRelation} and our assumptions that
$\La^2\equiv 0\pmod 2$ and $\delta\equiv i+j\pmod 2$ yield
$$
b_{i,j,k}(-2x)
=
(-1)^{\delta+i}
b_{i,j,k}(2x)
=
(-1)^j b_{i,j,k}(2x).
$$
Because $\La-\tilde w$ is characteristic, we have $(\La-\tilde w)^2\equiv \si\pmod 8$
and $\La^2\equiv \La\cdot(\La-\tilde w)\pmod 2$, so $\La\cdot \tilde w\equiv 0\pmod 2$.
Thus, $(\La-\tilde w)^2\equiv \si\pmod 8$ implies that $\La^2+\tilde w^2\equiv \si\pmod 4$
and so $\thalf(\tilde w^2-\si)\equiv \thalf\La^2\pmod 2$.
Therefore, by the definitions of $\eps(\tilde w,K_0)$ and $\tilde \eps(\tilde w,\La,K_0)$, we have
$$
\tilde \eps(\tilde w,\La,K_0)
-
\eps(\tilde w,K_0)
=
\frac{1}{2}(\tilde w^2-\si)-\frac{1}{2}K_0\cdot\La
\equiv
\frac{1}{2} (\La^2+K_0\cdot\La)
\pmod 2.
$$
By the preceding analysis, we can rewrite the coefficient
\eqref{eq:RHSCoeff2a} as
\begin{equation}
\label{eq:RHSCoeff2b}
\begin{aligned}
{}&
(-1)^{\eps(\tilde w,K_0)+x+y}
\binom{i}{i_0,\,\ldots,i_n}\SW'_X(K)
b_{i,j,k}(2x)
\\
&\quad\times
p^1(i_2,\dots,i_n)
\left(
(-1)^j-(-1)^{i_1}
\right).
\end{aligned}
\end{equation}
Lemma \ref{lem:AlgCoeff}
implies that the coefficients \eqref{eq:LHSCoefficientForBigi}
and \eqref{eq:RHSCoeff2b} must be equal.
For this to be a non-trivial
relation, we must have that
$p^1(i_2,\dots,i_n)$ is non-zero and consequently
we must have $i_u\equiv 1\pmod 2$ for $u=2,\dots,n$.
For $j$ even, take $i_1=\dots=i_n=1$ and $i_0=i-n$ while for $j$ odd, we take
$i_1=0$, $i_2=\dots=i_n=1$, and $i_0=i-n+1$ to get the desired equalities.
\end{proof}

\begin{rem}
\label{rmk:LowerCoeffProb}
Proposition \ref{prop:HighDegreeCoefficients} only determines the coefficients
$b_{i,j,k}(\chi_h,c_1^2,K\cdot\La,\La^2,m)$ for $i\ge \chi_h-c_1^2-3$.
An early manuscript version \cite{FLWConjecture_gap} of this article failed to note
that because $p^1(i_2,\dots,i_n)$ vanishes for low values of $i$ (since $i=i_0+i_1+\dots+i_n$ and so $i$ small implies that each $i_q$ is small)
the resulting relations were trivial and gave no information about the coefficients $b_{i,j,k}$.
\end{rem}

\begin{rem}[Determining the remaining coefficients]
We now describe some limitations on the ability of Equation
\eqref{eq:BlownUpUsefulCobordismFormula} to determine
the coefficients $b_{i,j,k}$ using the four-manifolds, $X_h$, constructed in Lemma \ref{lem:UsefulWithc=3}.
For $\chi_h$, $c_1^2$, $\La^2$, and $m$ fixed,
define a function $c_{i,j,k}:\ZZ\to \RR$ by setting
$c_{i,j,k}(x):=b_{i,j,k}(\chi_h,c_1^2,x,\La^2,m)$.
If, in the notation of Proposition \ref{prop:HighDegreeCoefficients}, one takes
$$
\La=y f_1 +f_2 +\sum_{u=1}^n \la_u e_u^*,
$$
then Lemma \ref{lem:PermutationSumAsDifferenceOperator} implies that
the coefficient of the term \eqref{eq:GeneralTerm} on the right-hand side
of \eqref{eq:BlownUpUsefulCobordismFormula} would be
$$
\nabla^{i_1+1}_{2\la_1}\dots \nabla^{i_n+1}_{2\la_n} c_{i,j,k}(K_0\cdot\La).
$$
Because $\nabla^{1}_{2\la_1}\cdots \nabla^{1}_{2\la_n} p(x)=0$
for any polynomial $p(x)$ of degree $n-1$ or less, the arguments used in the proof of 
Proposition \ref{prop:HighDegreeCoefficients} using the four-manifolds $X_h$
cannot determine the coefficients
$b_{0,j,k}$.  Arguing by induction on $i=v$ and 
by varying $i_1,\dots,i_v$,
one can show that the arguments of
Proposition \ref{prop:HighDegreeCoefficients} using the four-manifolds $X_h$ determine $b_{i,j,k}$ only up to a polynomial of degree $n-i-1$ in $\La\cdot K$.

This failure of
Proposition \ref{prop:HighDegreeCoefficients} to determine the coefficients $b_{i,j,k}$ using blow-ups of the manifolds $X_h$ stems from the failure of the set
$B'(\widetilde X_h(n))$ to be linearly independent.  Further progress with our method
would appear to rely on finding four-manifolds, $Y$, with $c(Y)>3$ and $B'(Y)$ admitting
few linear relations.  The `superconformal simple-type bound',
$$
c_1^2(Y)\ge \chi_h(Y) -2|B(Y)/\{\pm 1\}| -1,
$$
appearing in \cite[Theorem 4.1]{MMPdg} holds
for all known standard four-manifolds and indicates that the number of basic classes increases
as $c(Y)$ increases.  Consequently, one would need to search for standard four-manifolds where
the dimension of the span of $B'(Y)$ is large.
\end{rem}

\begin{proof}[Proof of Theorem \ref{thm:WittenSimpleType} for four-manifolds with $c_1^2\ge \chi_h -3$]
Assume that $Y$ is a standard four-manifold with $c_1^2(Y)\ge \chi_h(Y)-3$.
Let $X_h$ be a useful four-manifold provided by Lemma \ref{lem:UsefulWithc=3} with
$\chi_h(X_h)=\chi_h(Y)$.
By Theorem \ref{thm:WCBlowDownInvariance} and by blowing-up $Y$ if necessary,
we can assume that $c_1^2(Y)=c_1^2(X_h)$.  Let $\widetilde Y$ and $\widetilde X_h$ be the blow-ups of
$Y$ and $X_h$, respectively, at a point. Let $e^*\in H^2(\widetilde Y;\ZZ)$ be the Poincar\'e dual
of the exceptional curve.
For a  characteristic class $w\in H^2(Y;\ZZ)$, define
$\tilde w=w+e^*\in H^2(\widetilde Y;\ZZ)$.
Denoting $B'(Y)=\{K_1,\dots,K_b\}$,
there are cohomology classes $f_1,f_2\in H^2(Y;\ZZ)$ with $K_1\cdot f_i=0$ and
$f_i^2=0$ for $i=1,2$, and $f_1\cdot f_2=1$ by \cite[Corollary A.3]{FL2a}.
For a given $\delta$, we can choose an integer $a$ such that,
for $\La=2(af_1+f_2)\in H^2(Y;\ZZ)\subset H^2(\widetilde Y;\ZZ)$, we have $\La^2=8a$
and $I(\La)>\delta$.
Because $I(\La)>\delta$ and $\La-\tilde w$ is characteristic, we can use
this $\tilde w$ and $\La$ in
Lemma \ref{lem:ReduceCobordismFormToB'Sum}
to compute the degree-$\delta$ Donaldson invariant of $Y$.
Since $\La^2\equiv 0\pmod 2$ and $K_i$ is characteristic,
$K_i\cdot\La\equiv 0\pmod 2$ for all $K_i\in B(\widetilde Y)$.
Proposition \ref{prop:HighDegreeCoefficients} then only gives an expression for the
coefficients
$$
b_{i,j,k}((K_i\pm e^*)\cdot \La)
=
b_{i,j,k}(\chi_h(\tilde Y),c_1^2(\tilde Y),(K_i\pm e^*)\cdot\La,8a,m)
$$
appearing in Equation \eqref{eq:CompareCoeff2}
for $i\ge 1$.  We next show that we can ignore the terms
in Equation \eqref{eq:CompareCoeff2} with $i=0$.

As $\tilde w-\La$ is characteristic, we have
\begin{align*}
\tilde\eps(\tilde w,\La,K_i+e^*)
{}&\equiv
\tilde\eps(\tilde w,\La,K_i-e^*)+(\tilde w-\La)\cdot e^*\pmod 2\\
{}&\equiv
\tilde\eps(\tilde w,\La,K_i-e^*)+1\pmod 2.
\end{align*}
Using the fact that $(K_i+e^*)\cdot\La=(K_i-e^*)\cdot\La$, we obtain
$$
b_{i,j,k}((K_i+e^*)\cdot \La)=b_{i,j,k}((K_i-e^*)\cdot\La).
$$
Finally, because $n(K_i\pm e^*)=1$,
the terms for $K_i+e^*$ and $K_i-e^*$ in Equation \eqref{eq:CompareCoeff2} with $i=0$ will cancel out.
Thus, we may ignore the $i=0$ terms.

Since $\widetilde w$ is characteristic, the definition of $\tilde \eps$
in \eqref{eq:DefineTildeEps} implies that
$$
{\tilde\eps(\tilde w,\La,K_i\pm e^*)}+\frac{1}{2}\La\cdot (K_i\pm e^*)
\equiv
\eps(\tilde w,K_i\pm e^*)
\pmod 2.
$$
Therefore, the formula for the coefficients, $b_{i,j,k}$, in
Proposition \ref{prop:HighDegreeCoefficients} and the vanishing of the
terms with $i=0$ allow us to rewrite Equation \eqref{eq:CompareCoeff2} as
\begin{equation}
\begin{aligned}
\label{eq:CompareCoeff2a}
{}&
D^{\tilde w}_{\widetilde Y}(h^{\delta-2m}x^m)
\\
{}&\quad=
\sum_{\begin{subarray}{l}i+j+2k\\=\delta-2m\end{subarray}}
\sum_{K\in B'(\widetilde Y)}
(-1)^{{\eps(\tilde w,K)}}SW'_X(K)
\frac{(\delta-2m)!}{k!i!}2^{m-k-1}
\langle K,h\rangle^i
Q_X^k(h).
\end{aligned}
\end{equation}
Comparing Equations \eqref{eq:CompareCoeff2a} and
\eqref{eq:DInvarForWCB'Sum}, noting that $c(\widetilde Y)=4$, and applying Lemma \ref{lem:ReduceDFormToB'Sum}
then shows that Witten's Conjecture \ref{conj:WittenSimpleType} holds for $\widetilde Y$ and thus for $Y$.
\end{proof}

Before proceeding to the proof of Theorem \ref{thm:WittenSimpleType} for abundant four-manifolds, we recall 
a vanishing result for abundant four-manifolds. If $Y$ is a standard four-manifold, $w\in H^2(Y;\ZZ)$, and $h\in H_2(Y;\RR)$, we define
$$
\SW_{Y,i}^w(h)
:=
\sum_{K\in B(Y)} (-1)^{\eps(w,K)}\SW'_Y(K)\langle K,h\rangle^i.
$$
We then
recall the

\begin{thm}
\label{thm:LowDegVanishingForAbundant}
\cite[Theorem 1.1]{FKLM}
Theorem \ref{thm:Cobordism} implies
that if $Y$ is a standard and abundant
four-manifold and $w$ is characteristic,
then $\SW^w_{Y,i}$ vanishes for $i<c(Y)-2$.
\end{thm}

\begin{proof}[Proof of Theorem \ref{thm:WittenSimpleType} for abundant four-manifolds]
We now show that Proposition \ref{prop:HighDegreeCoefficients} suffices to
prove Witten's Conjecture \ref{conj:WittenSimpleType} for abundant four-manifolds.
By the argument in the proof of Lemma \ref{lem:ReduceDFormToB'Sum},
for $w$ characteristic (so $w^2\equiv c_1^2(Y)\pmod 2$),
\begin{equation}
\label{eq:SWPolyB'}
\SW_{Y,i}^w(h)
=
(1+(-1)^{c(X)+i})\sum_{K\in B'(Y)}
(-1)^{\eps(w,K)}n(K)\SW'_Y(K)\langle K,h\rangle^i.
\end{equation}
By Theorem \ref{thm:WCBlowDownInvariance}, it suffices to prove that
Conjecture \ref{conj:WittenSimpleType} holds for  the blow-up of $Y$
at any number of points.  We can therefore assume that $c_1^2(Y)=\chi_h(Y)-3-n$ for $n\ge 1$.
For any non-negative integers $\delta$ and $m$ satisfying
$\delta-2m\ge 0$, choose an integer $a$ such that $8a>\delta-5\chi_h(Y)-c_1^2(Y)$.
Let $f_1,f_2\in B(Y)^\perp$ satisfy $f_1\cdot f_2=1$ and $f_i^2=0$.
Then for $\La=2af_1+2f_2$, we have $I(\La)>\delta$ as required
in Lemma \ref{lem:ReduceCobordismFormToB'Sum}.
Note that because $\La\equiv 0\pmod 2$, for $w$ characteristic, the class $\La-w$
is also characteristic.  Since $w$ is characteristic and $\La\in B(Y)^\perp$,
we have
$$
\tilde \eps(w,\La,K)\equiv\eps(w,K) \pmod 2.
$$
For $\La\in B'(Y)^\perp$, we have $b_{i,j,k}=0$ unless $c(Y)+i\equiv 0\pmod 2$
by \eqref{eq:CombinedCoeffReverseLaKEquality} and hence
$1+(-1)^{c(X)+i}=2$ in Equation \eqref{eq:SWPolyB'}.
As $K\cdot \La$ and hence the coefficients $b_{i,j,k}=b_{i,j,k}(K\cdot\La)$
are independent of
$K\in B'(Y)$, we can write the expression for the Donaldson invariant
in Lemma \ref{lem:ReduceCobordismFormToB'Sum} as
\begin{equation}
\label{eq:CompareCoeffAbund}
D^w_Y(h^{\delta-2m}x^m)=
\sum_{\begin{subarray}{l}i+j+2k\\=\delta-2m\end{subarray}}
\frac{1}{2} b_{i,j,k} \SW^w_{Y,i}(h)
\langle \La,h\rangle^j
Q_Y(h)^k.
\end{equation}
Theorem \ref{thm:LowDegVanishingForAbundant} allows us to ignore the coefficients
$b_{i,j,k}$ in Equation \eqref{eq:CompareCoeffAbund} with $i\le n=c(Y)-3$.
By Proposition \ref{prop:HighDegreeCoefficients},
we then can
rewrite Equation \eqref{eq:CompareCoeffAbund} as
$$
D^w_Y(h^{\delta-2m}x^m)
=
\sum_{\begin{subarray}{l}i+2k\\=\delta-2m\end{subarray}}
\sum_{K\in B'(Y)}
(-1)^{\eps(w,K)}
\frac{(\delta-2m)!}{2^{n+k-m}k!i!} n(K)\SW'_Y(K)\langle K,h\rangle^i
Q_Y(h)^k.
$$
Comparing this expression for $D^w_Y(h^{\delta-2m}x^m)$ with that in 
Equation \eqref{eq:DInvarForWCB'Sum}
then completes the proof of the
theorem.
\end{proof}

\bigskip
\footnotesize
\noindent\textit{Acknowledgments.}
Feehan was supported in part by NSF grant DMS 0125170.  Leness was supported in part by a Florida International University Summer Research Grant and NSF grant DMS 0905786. Leness is indebted to Ron Stern and Ron Fintushel for considerable help with examples, to Nick Saveliev for comments on drafts, and to Miroslav Yotov for Lemma \ref{lem:AlgCoeff}. In addition, Leness would like to thank the organizers of the Park City Mathematics Institute 2006 Summer School for providing an outstanding research environment. Feehan is grateful to Brendan Owens for help with questions on knot theory. Both authors warmly thank Tom Mrowka for his faithful encouragement of this project since its inception in 1994 and are very grateful to Yasha Eliashberg for his encouragement and steadfast support while we prepared the final version of this article.

%
%

\def\cprime{$'$} \def\polhk#1{\setbox0=\hbox{#1}{\ooalign{\hidewidth
  \lower1.5ex\hbox{`}\hidewidth\crcr\unhbox0}}} \def\cprime{$'$}
  \def\cprime{$'$} \def\cprime{$'$}
  \def\lfhook#1{\setbox0=\hbox{#1}{\ooalign{\hidewidth
  \lower1.5ex\hbox{'}\hidewidth\crcr\unhbox0}}} \def\cprime{$'$}
  \def\cprime{$'$} \def\cprime{$'$} \def\cprime{$'$}
\providecommand{\bysame}{\leavevmode\hbox to3em{\hrulefill}\thinspace}
\providecommand{\MR}{\relax\ifhmode\unskip\space\fi MR }
\providecommand{\MRhref}[2]{%
  \href{http://www.ams.org/mathscinet-getitem?mr=#1}{#2}
}
\providecommand{\href}[2]{#2}

\begin{comment}
\bibliography{master,mfpde}
\bibliographystyle{amsplain}

\end{document}